\documentclass{amsart}

\usepackage{amsfonts, amssymb, amsmath, stmaryrd}
\usepackage{mathrsfs,array}
\usepackage[pdftex]{graphicx}
\usepackage{tikz}
\usetikzlibrary{fadings}
\usepackage{xy}
\usepackage{hyperref}
\usepackage{verbatim}
\input xy
\xyoption{all}

\newtheorem{theorem}{Theorem}[section]
\newtheorem{corollary}[theorem]{Corollary}

\newtheorem{proposition}[theorem]{Proposition}
\newtheorem{conjecture}[theorem]{Conjecture}
\newtheorem{construction}[theorem]{Construction}

\newtheorem{definition}[theorem]{Definition}
\theoremstyle{definition}
\newtheorem{remark}[theorem]{Remark}

\DeclareMathOperator{\Gr}{Gr}
\DeclareMathOperator{\Spa}{Spa}
\DeclareMathOperator{\Spd}{Spd}
\DeclareMathOperator{\Spec}{Spec}
\DeclareMathOperator{\Spf}{Spf}
\DeclareMathOperator{\Fil}{Fil}

\DeclareMathOperator{\Gal}{Gal}
\DeclareMathOperator{\Res}{Res}

\DeclareMathOperator{\Lie}{Lie}

\DeclareMathOperator{\Kt}{Kt}
\DeclareMathOperator{\FF}{FF}

\def\et{\mathrm{\acute{e}t}}

\def\inf{\mathrm{inf}}

\def\dR{\mathrm{dR}}

\def\perf{\mathrm{perf}}

\renewcommand{\diamond}{\diamondsuit}

\newcommand{\Z}{\mathbb{Z}}
\newcommand{\F}{\mathbb{F}}
\newcommand{\Q}{\mathbb{Q}}
\newcommand{\R}{\mathbb{R}}
\newcommand{\A}{\mathbb{A}}

\newcommand{\C}{\mathbb{C}}
\newcommand{\G}{\mathbb{G}}
\renewcommand{\P}{\mathbb{P}}
\newcommand{\OO}{\mathcal{O}}

\newcommand{\GL}{\mathrm{GL}}

\newcommand{\End}{\mathrm{End}}

\newcommand{\THH}{\mathrm{THH}}
\newcommand{\TC}{\mathrm{TC}}

\newcommand{\Fl}{\mathscr{F}\!\ell}

\newcommand{\Perf}{\mathrm{Perf}}
\newcommand{\Sh}{\mathrm{Sh}}
\newcommand{\Frob}{\mathrm{Frob}}

\title{$p$-adic geometry}

\author{Peter Scholze}
\address{Mathematisches Institut der Universit\"at Bonn\\ Endenicher Allee 60\\ 53115 Bonn\\ Germany}
\email{scholze@math.uni-bonn.de}

\begin{document}

\begin{abstract}
We discuss recent developments in $p$-adic geometry, ranging from foundational results such as the degeneration of the Hodge-to-de Rham spectral sequence for ``compact $p$-adic manifolds'' over new period maps on moduli spaces of abelian varieties to applications to the local and global Langlands conjectures, and the construction of ``universal'' $p$-adic cohomology theories. We finish with some speculations on how a theory that combines all primes $p$, including the archimedean prime, might look like.
\end{abstract}



\maketitle

\section{Introduction}

In this survey paper, we want to give an introduction to the world of ideas which the author has explored in the 
past few years, and indicate some possible future directions. The two general themes that dominate this work are the cohomology of algebraic varieties, and the local and global Langlands correspondences. These two topics are classically intertwined ever since the cohomology of the moduli space of elliptic curves and more general Shimura varieties has been used for the construction of Langlands correspondences. Most of our work so far is over $p$-adic fields, where we have established analogues of the basic results of Hodge theory for ``compact $p$-adic manifolds'', have constructed a ``universal'' $p$-adic cohomology theory, and have made progress towards establishing the local Langlands correspondence for a general $p$-adic reductive group by using a theory of $p$-adic shtukas, and we will recall these results below.

However, here we wish to relay another, deeper, relation between the cohomology of algebraic varieties and the structures underlying the Langlands corresondence, a relation that pertains not to the cohomology of specific algebraic varieties, but to the very notion of what ``the'' cohomology of an algebraic variety is. Classically, the study of the latter is the paradigm of ``motives'' envisioned by Grothendieck; however, that vision has still only been partially realized, by Voevodsky, \cite{VoevodskyMotives}, and others. Basically, Grothendieck's idea was to find the ``universal'' cohomology as the universal solution to a few basic axioms; in order to see that this has the desired properties, one however needs to know the existence of ``enough'' algebraic cycles as encoded in the standard conjectures, and more generally the Hodge and Tate conjectures. However, little progress has been made on these questions. We propose to approach the subject from the other side and construct an \emph{explicit} cohomology theory that practically behaves like a universal cohomology theory (so that, for example, it specializes to all other known cohomology theories); whether or not it is universal in the technical sense of being the universal solution to certain axioms will then be a secondary question.

This deeper relation builds on the realization of Drinfeld, \cite{DrinfeldGL2}, that in the function field case, at the heart of the Langlands correspondence lie moduli spaces of \emph{shtukas}. Anderson, \cite{Anderson}, Goss, \cite{Goss}, and others have since studied the notion of t-motives, which is a special kind of shtuka, and is a remarkable function field analogue of motives, however without any relation to the cohomology of algebraic varieties. What we are proposing here is that, despite extreme difficulties in making sense of this, there should exist a theory of shtukas in the number field case, and that the cohomology of an algebraic variety, i.e.~a motive, should be an example of such a shtuka.

This picture has been essentially fully realized in the $p$-adic case. In the first sections of this survey, we will explain these results in the $p$-adic case; towards the end, we will then speculate on how the full picture over $\Spec \Z$ should look like, and give some evidence that this is a reasonable picture.

{\bf Acknowledgments.} This survey was written in relation to the author's lecture at the ICM 2018. Over the past years, I have benefitted tremendously from discussions with many mathematicians, including Bhargav Bhatt, Ana Caraiani, Gerd Faltings, Laurent Fargues, Ofer Gabber, Eugen Hellmann, Lars Hesselholt, Kiran Kedlaya, Mark Kisin, Arthur-C\'esar le Bras, Akhil Mathew, Matthew Morrow, Wies\l{}awa Nizio\l{}, Michael Rapoport, Richard Taylor, and many others. It is a great pleasure to thank all of them for sharing their insights, and more generally the whole mathematical community for being so generous, inviting, supportive, challenging, and curious. I feel very glad and honoured to be part of this adventure.

\section{$p$-adic Hodge theory}

Let us fix a complete algebraically closed extension $C$ of $\Q_p$. The analogue of a compact complex manifold in this setting is a proper smooth rigid-analytic variety $X$ over $C$.\footnote{Rigid-analytic varieties were first defined by Tate, \cite{TateRigid}, and alternative and more general foundations have been proposed by various authors, including Raynaud, \cite{Raynaud}, Berkovich, \cite{Berkovich}, Fujiwara-Kato, \cite{FujiwaraKato}, and Huber, \cite{HuberBook}. We have found Huber's setup to be the most natural, and consequently we will often implicitly regard all schemes, formal schemes and rigid spaces that appear in the following as adic spaces in the sense of Huber; his category naturally contains all of these categories as full subcategories.} Basic examples include the analytification of proper smooth algebraic varieties over $C$, or the generic fibres of proper smooth formal schemes over the ring of integers $\OO_C$ of $C$. More exotic examples with no direct relation to algebraic varieties are given by the Hopf surfaces
\[
X= ((\mathbb A^2_C)^{\mathrm{rig}}\setminus \{(0,0)\})/q^\Z\ ,
\]
where $q\in C$ is an element with $0<|q|<1$ acting via diagonal multiplication on $\mathbb A^2$. Recall that their complex analogues are non-K\"ahler, as follows from the asymmetry in Hodge numbers, $H^0(X,\Omega_X^1)=0$ while $H^1(X,\OO_X)=C$. On the other hand, some other non-K\"ahler manifolds such as the Iwasawa manifolds do not have a $p$-adic analogue.

The basic cohomological invariants of a compact complex manifold also exist in this setting. The analogue of singular cohomology is \'etale cohomology $H^i_\et(X,\Z_\ell)$, which is defined for any prime $\ell$, including $\ell=p$. If $\ell\neq p$, it follows for example from work of Huber, \cite[Proposition 0.5.3]{HuberBook}, that this is a finitely generated $\Z_\ell$-module. For $\ell=p$, this is also true by \cite[Theorem 1.1]{ScholzePAdicHodge}, but the argument is significantly harder.

Moreover, one has de~Rham cohomology groups $H^i_\dR(X/C)$ and Hodge cohomology groups $H^i(X,\Omega_X^j)$, exactly as for compact complex manifolds. These are finite-dimensional by a theorem of Kiehl, \cite{KiehlProper}. By the definition of $H^i_\dR(X/C)$ as the hypercohomology of the de~Rham complex, one finds an $E_1$-spectral sequence
\[
E_1^{ij} = H^j(X,\Omega_X^i)\Rightarrow H^{i+j}_\dR(X/C)
\]
called the Hodge-to-de~Rham spectral sequence. In complex geometry, a basic consequence of Hodge theory is that this spectral sequence degenerates at $E_1$ if $X$ admits a K\"ahler metric. This assumption is not necessary in $p$-adic geometry:

\begin{theorem}[{\cite[Corollary 1.8]{ScholzePAdicHodge}}, {\cite[Theorem 13.12]{BMS}}] For any proper smooth rigid-analytic space $X$ over $C$, the Hodge-to-de~Rham spectral sequence
\[
E_1^{ij} = H^j(X,\Omega_X^i)\Rightarrow H^{i+j}_\dR(X/C)
\]
degenerates at $E_1$. Moreover, for all $i\geq 0$,
\[
\sum_{j=0}^i \dim_C H^{i-j}(X,\Omega_X^j) = \dim_C H^i_\dR(X/C) = \dim_{\Q_p} H^i_\et(X,\Q_p)\ .
\]
\end{theorem}

Fortunately, the Hodge-to-de~Rham spectral sequence does degenerate for the Hopf surface -- and the examples of nondegeneration such as the Iwasawa manifolds do not have $p$-adic analogues.

Over the complex numbers, the analogue of the equality $\dim_C H^i_\dR(X/C) = \dim_{\Q_p} H^i_\et(X,\Q_p)$ follows from the comparison isomorphism between singular and de~Rham cohomology. In the $p$-adic case, the situation is slightly more complicated, and the comparison isomorphism only exists after extending scalars to Fontaine's field of $p$-adic periods $B_\dR$. If $X$ is only defined over $C$, it is nontrivial to formulate the correct statement, as there is no natural map $C\to B_\dR$ along which one can extend scalars; the correct statement is Theorem~\ref{thm:BdRpluscohom} below.

There is however a different way to obtain the desired equality of dimensions. This relies on the Hodge-Tate spectral sequence, a form of which is implicit in Faltings's proof of the Hodge-Tate decomposition, \cite{FaltingsJAMS}.

\begin{theorem}[{\cite[Theorem 3.20]{ScholzeCDM}}, {\cite[Theorem 13.12]{BMS}}] For any proper smooth rigid-analytic space $X$ over $C$, there is a Hodge-Tate spectral sequence
\[
E_2^{ij} = H^i(X,\Omega_X^j)(-j)\Rightarrow H^{i+j}_\et(X,\Z_p)\otimes_{\Z_p} C
\]
that degenerates at $E_2$.
\end{theorem}

Here, $(-j)$ denotes a Tate twist, which becomes important when one wants to make everything Galois-equivariant. Note that the Hodge cohomology groups appear in the other order than in the Hodge-to-de~Rham spectral sequence.

\begin{remark} If $X$ is the base change of a proper smooth rigid space defined over a discretely valued field $K\subset C$, then everything in sight carries a Galois action, and it follows from the results of Tate, \cite{TatePDivGroups}, that there is a unique Galois-equivariant splitting of the abutment filtration, leading to a Galois-equivariant isomorphism
\[
H^i_\et(X,\Q_p)\otimes_{\Q_p} C\cong \bigoplus_{j=0}^i H^{i-j}(X,\Omega_X^j)(-j)\ ,
\]
answering a question of Tate, \cite[Section 4.1, Remark]{TatePDivGroups}. However, this isomorphism does not exist in families. This is analogous to the Hodge decomposition over the complex numbers that does not vary holomorphically in families.
\end{remark}

An interesting question is whether Hodge symmetry could still hold under some condition on $X$. Such an analogue of the K\"ahler condition has recently been proposed by Li, \cite{Li}. In joint work with Hansen, \cite{HansenLi}, they state the following conjecture that they prove in the case $i+j=1$.

\begin{conjecture} Let $X$ be a proper smooth rigid-analytic variety that admits a formal model $\mathfrak X$ whose special fibre is projective. Then for all $i,j\geq 0$, one has $\dim_C H^i(X,\Omega_X^j) = \dim_C H^j(X,\Omega_X^i)$.
\end{conjecture}

The condition is indeed analogous to the K\"ahler condition; in Arakelov theory, the analogue of a metric is a formal model, and the positivity condition on the K\"ahler metric finds a reasonable analogue in the condition that the special fibre is projective.

Note that in particular it follows from the results of Hansen-Li that for any proper formal model $\mathfrak X$ of the Hopf surface (which exist by Raynaud's theory of formal models), the special fibre is a non-projective (singular) proper surface.

From the Hodge-to-de~Rham and the Hodge-Tate spectral sequence, one obtains abutments filtrations that we call the Hodge-de~Rham filtration and the Hodge-Tate filtration. Their variation in families defines interesting period maps as we will recall in the next sections.

\section{Period maps from de~Rham cohomology}

First, we recall the more classical case of the period maps arising from the variation of the Hodge filtration on de~Rham cohomology. For simplicity, we will discuss the moduli space $\mathcal M/\mathbb Z$ of elliptic curves. If $E$ is a an elliptic curve over the complex numbers, then $H^0(E,\Omega_E^1) = (\Lie E)^\ast$ is the dual of the Lie algebra, $H^1(E,\OO_E)=\Lie E^\ast$ is the Lie algebra of the dual elliptic curve $E^\ast$ (which, for elliptic curves, is canonically isomorphic to $E$ itself), and the Hodge-de~Rham filtration is a short exact sequence
\[
0\to (\Lie E)^\ast\to H^1_\dR(E/\C)\to \Lie E^\ast\to 0\ ,
\]
where $H^1_\dR(E/\C)\cong H^1_{\mathrm{sing}}(E,\Z)\otimes_{\Z} \C$. The classical period map takes the form
\[\xymatrixrowsep{0.7pc}\xymatrix{
& \widetilde{\mathcal M}\ar[dl]\ar[dr]\\
\mathcal M(\C) && \mathbb H^{\pm} = \mathbb P^1(\C)\setminus \mathbb P^1(\mathbb R)
}\]
where $\widetilde{\mathcal M}\to \mathcal M(\C)$ is the $\GL_2(\Z)$-torsor parametrizing trivializations of the first singular cohomology of the elliptic curve. Given an elliptic curve $E/\C$ with such a trivialization $H^1(E,\Z)\cong \Z^2$, the Hodge filtration $(\Lie E)^\ast\subset H^1_{\mathrm{sing}}(E,\Z)\otimes_{\Z} \C\cong \C^2$ defines a point of $\mathbb P^1(\C)\setminus \mathbb P^1(\R)$ as
\[
H^1_{\mathrm{sing}}(E,\Z)\otimes_\Z \C\cong H^1_\dR(E/\C)\cong (\Lie E)^\ast\oplus \overline{(\Lie E)^\ast}\ .
\]

If now $E$ is an elliptic curve over the algebraically closed $p$-adic field $C$, then we still have the identifications $H^0(E,\Omega_E^1) = (\Lie E)^\ast$, $H^1(E,\OO_E)=\Lie E^\ast$ and the Hodge-de~Rham filtration
\[
0\to (\Lie E)^\ast\to H^1_\dR(E/C)\to \Lie E^\ast\to 0\ .
\]
To construct period maps, we need to fix a point $x\in \mathcal M(\overline{\F}_p)$. Let $\breve{\Q}_p$ be the completion of the maximal unramified extension of $\Q_p$, and consider the rigid space
\[
\breve{\mathcal M} = \mathcal M_{\breve{\Q}_p}^{\mathrm{rig}}\ .
\]
Let $U_x\subset \breve{\mathcal M}$ be the tube of $x$, i.e.~the open subspace of all points specializing to $x$. This is isomorphic to an open unit ball $\mathbb D = \{z\mid |z|<1\}$. Then the crystalline nature of de~Rham cohomology (in particular, the connection) imply that the de~Rham cohomology $H^i_\dR(E_{\tilde{x}})$ is canonically identified for all $\tilde{x}\in U_x$. Now the variation of the Hodge filtration defines a period map
\[
\pi_x: U_x\cong \mathbb D\to \breve{\mathbb P}^1 = (\mathbb P^1_{\breve{\Q}_p})^{\mathrm{rig}}\ .
\]

There are two cases to consider here. If $x$ corresponds to an ordinary elliptic curve, then there is an identification $U_x\cong \mathbb D$ such that the map $\pi_x$ is given by the logarithm map
\[
z\mapsto \log(1+z): U_x=\mathbb D = \{z\mid |z|<1\}\to \breve{\mathbb A}^1\subset \breve{\mathbb P}^1\ .
\]
This map is an \'etale covering map onto $\breve{\mathbb A}^1$, and the geometric fibres are given by copies of $\Q_p/\Z_p$; for example, the fibre over $0$ is given by all $p$-power roots of unity minus $1$, i.e.~$z=\zeta_{p^r}-1$. Note that this implies in particular that $(\mathbb A^1_{\C_p})^{\mathrm{rig}}$ has interesting (non-finite) \'etale coverings, contrary both to the scheme case and the case over the complex numbers.

If $x$ corresponds to a supersingular elliptic curve, then $\pi_x$ becomes a map
\[
\pi_x: U_x\cong \mathbb D\to \breve{\mathbb P}^1\ .
\]
In this case, the map is known as the Gross-Hopkins period map, \cite{GrossHopkins}. It is an \'etale covering map of $\breve{\mathbb P}^1$ whose geometric fibres are given by copies of $\GL_2(\Q_p)_1/\GL_2(\Z_p)$, where $\GL_2(\Q_p)_1\subset \GL_2(\Q_p)$ is the open subgroup of all $g\in \GL_2(\Q_p)$ with $\det g\in \Z_p^\times$. 

It is important to note that these maps $\pi_x$ for varying $x$ cannot be assembled into a single map from $\breve{\mathcal M}$ towards $\breve{\mathbb P}^1$, contrary to the case over the complex numbers. We will however see a global period map in the next section.

These examples are the basic examples of \emph{local Shimura varieties}, which are associated with local Shimura data, \cite{RapoportViehmann}:
\begin{enumerate}
\item[{\rm (i)}] A reductive group $G$ over $\Q_p$.
\item[{\rm (ii)}] A conjugacy class of minuscule cocharacters $\mu: \G_m\to G_{\overline{\Q}_p}$, defined over the reflex field $E$ (a finite extension of $\Q_p$).
\item[{\rm (iii)}] A $\sigma$-conjugacy class $b\in B(G,\mu)$.
\end{enumerate}

We will say something about datum (iii) in Theorem~\ref{thm:BGCurve} below. In our example, we have $G=\GL_2$, and $\mu$ is the conjugacy class of $t\mapsto \mathrm{diag}(t,1)$. In this case, $B(G,\mu)$ contains exactly two elements, corresponding to the cases of ordinary and supersingular elliptic curves, respectively.

Let $\breve{E}$ be the completion of the maximal unramified extension of $E$. Corresponding to $G$ and $\mu$, one gets a flag variety $\breve{\Fl}_{G,\mu}$ over $\breve{E}$, which we consider as an adic space; in our example, this is $\breve{\mathbb P}^1$. The following theorem on the existence of local Shimura varieties proves a conjecture of Rapoport-Viehmann, \cite{RapoportViehmann}.

\begin{theorem}[\cite{Berkeley}]\label{thm:localshimura} There is a natural open subspace $\Fl_{G,\mu}^a\subset \Fl_{G,\mu}$, called the admissible locus, and an \'etale covering map
\[
\mathcal M_{(G,b,\mu),K}\to \Fl_{G,\mu}^a\subset \Fl_{G,\mu}
\]
with geometric fibres $G(\Q_p)/K$ for any compact open subgroup $K\subset G(\Q_p)$.
\end{theorem}

\begin{remark} The geometric fibres differ from the examples above. In the examples, $\mathcal M_{(G,b,\mu),\GL_2(\Z_p)}$ is a certain disjoint union of various $U_x$ for $x$ ranging through a $p$-power-isogeny class of elliptic curves.
\end{remark}

The (towers of) spaces $(\mathcal M_{(G,b,\mu),K})_{K\subset G(\Q_p)}$ are known as local Shimura varieties; their cohomology is expected to realize local Langlands correspondences, cf.~\cite{RapoportViehmann}. This is one of our primary motivations for constructing these spaces.

In hindsight, one can say that in \cite{RapoportZink}, Rapoport-Zink constructed these local Shimura varieties in many examples, by explicitly constructing $\mathcal M_{(G,b,\mu),K}$ as a moduli space of $p$-divisible groups with extra structures such as endomorphisms and polarization, together with a quasi-isogeny to a fixed $p$-divisible group. An advantage of this approach is that, at least for special choices of $K$ such as parahoric subgroups, one actually constructs formal schemes whose generic fibre is the local Shimura variety, and it is often easier to understand the formal scheme.

We will explain our construction of local Shimura varieties in Section~\ref{sec:twistor}, and in Section~\ref{sec:padicshtuka} we will deal with integral models of local Shimura varieties.

\section{Period maps from \'etale cohomology}

A different period map known as the Hodge-Tate period map parametrizes the variation of the Hodge-Tate filtration in families, and has been defined for general Shimura varieties of Hodge type in \cite{ScholzeTorsion}.

If $E$ is an elliptic curve over the algebraically closed $p$-adic field $C$, the Hodge-Tate filtration is given by a short exact sequence
\[
0\to \Lie E^\ast\to H^1_\et(E,\Z_p)\otimes_{\Z_p} C\to (\Lie E)^\ast(-1)\to 0\ .
\]
Note that the Lie algebra terms here appear in opposite order when compared to the situation over $\C$. To obtain the period map, we need to trivialize the middle term. Thus, using the space $\breve{\mathcal M} = (\mathcal M_{\breve{\Q}_p})^{\mathrm{rig}}$ as in the last section, we consider the diagram
\[\xymatrixrowsep{0.7pc}\xymatrix{
& \widetilde{\breve{\mathcal M}}\ar[dl]\ar[dr]^{\pi_{HT}} \\
\breve{\mathcal M} && \breve{\mathbb P}^1\ .
}\]
Here $\widetilde{\breve{\mathcal M}}\to \breve{\mathcal M}$ parametrizes isomorphisms $H^1_\et(E,\Z_p)\cong \Z_p^2$; this defines a $\GL_2(\Z_p)$-torsor. An essential difficulty here is that this space will be very big, and in particular highly nonnoetherian. By \cite{ScholzeTorsion}, it gives a basic example of a perfectoid space in the sense of \cite{ScholzeThesis}.

Now $\widetilde{\breve{\mathcal M}}$ admits the Hodge-Tate period map
\[
\pi_{HT}: \widetilde{\breve{\mathcal M}}\to \breve{\mathbb P}^1
\]
sending a pair of an elliptic curve $E/C$ with an isomorphism $H^1_\et(E,\Z_p)\cong \Z_p^2$ to the filtration $\Lie E^\ast\subset H^1_\et(E,\Z_p)\otimes_{\Z_p} C\cong C^2$.

The geometry of this map is very interesting. When restricted to Drinfeld's upper half-plane $\breve{\Omega}^2 = \breve{\mathbb P}^1\setminus \mathbb P^1(\Q_p)$, the map is a pro-finite \'etale cover, while the fibres over points in $\mathbb P^1(\Q_p)$ are curves, so that the fibre dimension jumps. Let us discuss these two situations in turn.

If we fix a supersingular elliptic curve $x\in \mathcal M(\overline{\F}_p)$, then we can restrict the Hodge-Tate period map to
\[
\widetilde{U}_x = U_x\times_{\breve{\mathcal M}} \widetilde{\breve{\mathcal M}}\ ,
\]
and arrive at the following picture:\vspace{-0.1in}
\[\xymatrixrowsep{0pc}\xymatrix{
&& \widetilde{U}_x\ar[dl]\ar^{\pi_{HT}}[ddr]\\
& U_x\ar[dl]_{\pi_x} \\
\breve{\mathbb P}^1 &&& \breve{\Omega}^2\ .
}\]
Here, the map $\widetilde{U}_x\to U_x$ is a $\GL_2(\Z_p)$-torsor, while $U_x\to \breve{\mathbb P}^1$ has geometric fibers $\GL_2(\Q_p)_1/\GL_2(\Z_p)$; in fact, in total $\widetilde{U}_x\to \breve{\mathbb P}^1$ is a $\GL_2(\Q_p)_1$-torsor. To restore full $\GL_2(\Q_p)$-equivarience, we consider the perfectoid space
\[
\mathcal M_{LT,\infty} = \widetilde{U}_x\times^{\GL_2(\Q_p)_1} \GL_2(\Q_p)\cong \bigsqcup_\Z \widetilde{U}_x\ ,
\]
which is known as the Lubin-Tate tower at infinite level, cf.~\cite{ScholzeWeinstein}, \cite{WeinsteinLubinTate}. On the other hand, $\widetilde{U}_x\to \breve{\Omega}^2$ turns out to be an $\OO_D^\times$-torsor, where $D/\Q_p$ is the quaternion algebra. Here, the $\OO_D$-action arises from the identification $\OO_D = \mathrm{End}(E_x)\otimes_{\Z} \Z_p$, where $E_x$ is the supersingular elliptic curve corresponding to $x$; by functoriality, this acts on the deformation space of $x$, and thus also on $U_x$ and $\widetilde{U}_x$. In terms of the Lubin-Tate tower,
\[
\pi_{HT}: \mathcal M_{LT,\infty}\to \breve{\Omega}^2
\]
is a $D^\times$-torsor.

We see that $\mathcal M_{LT,\infty}$ has two different period morphisms, corresponding to the Hodge-de~Rham filtration and the Hodge-Tate filtration. This gives the isomorphism between the Lubin-Tate and Drinfeld tower at infinite level, \cite{FaltingsTwoTowers}, \cite{FarguesTwoTowers},
\[\xymatrixrowsep{0pc}\xymatrix{
& \mathcal M_{LT,\infty}\ar@{=}[r]\ar[dl]_-{\GL_2(\Q_p)} & \mathcal M_{Dr,\infty}\ar[dr]^{D^\times}\\
\breve{\mathbb P}^1 &&& \breve{\Omega}^2\ .
}\]

A similar duality theorem holds true for any local Shimura variety for which $b$ is basic. The local Shimura datum $(G,b,\mu)$ then has a dual datum $(\widehat{G},\widehat{b},\widehat{\mu})$, where $\widehat{G} = J_b$ is the $\sigma$-centralizer of $b$, $\widehat{b} = b^{-1}\in J_b$, and $\widehat{\mu} = \mu^{-1}$ under the identification $G_{\overline{\Q}_p}\cong \widehat{G}_{\overline{\Q}_p}$.

\begin{theorem}[\cite{ScholzeWeinstein}, \cite{Berkeley}] There is a natural isomorphism
\[
\varprojlim_{K\subset G(\Q_p)} \mathcal M_{(G,b,\mu),K}\cong \varprojlim_{\widehat{K}\subset\widehat{G}(\Q_p)} \mathcal M_{(\widehat{G},\widehat{b},\widehat{\mu}),\widehat{K}}\ .
\]
\end{theorem}

\begin{remark} This proves a conjecture of Rapoport-Zink, \cite[Section 5.54]{RapoportZink}. One has to be careful with the notion of inverse limits here, as inverse limits in adic spaces do not exist in general. One interpretation is to take the inverse limit in the category of diamonds discussed below.
\end{remark}

On the other hand, we can restrict the Hodge-Tate period map to $\mathbb P^1(\Q_p)$. In this case, the fibres of $\pi_{HT}$ are curves. More precisely, consider the Igusa curve $\mathrm{Ig}$ parametrizing ordinary elliptic curves $E$ in characteristic $p$ together with an isomorphism $E[p^\infty]_\et\cong \Q_p/\Z_p$ of the of the \'etale quotient $E[p^\infty]_\et$ of the $p$-divisible group $E[p^\infty]$. Its perfection $\mathrm{Ig}^\perf$ lifts uniquely to a flat formal scheme $\mathfrak{Ig}$ over $\breve{\Z}_p$, and then, modulo boundary issues, $\pi_{HT}^{-1}(x)$ is the generic fibre of $\mathfrak{Ig}$.

Let us state these results about the Hodge-Tate period map for a general Shimura variety of Hodge type.

\begin{theorem}[\cite{ScholzeTorsion}, \cite{CaraianiScholze}]\label{thm:fibreigusa} Consider a Shimura variety $\Sh_K$, $K\subset G(\A_f)$, of Hodge type, associated with some reductive group $G/\Q$ and Shimura data, including the conjugacy class of cocharacters $\mu$ with field of definition $E$. Let $\Fl_{G,\mu}=G/P_\mu$ be the corresponding flag variety.\footnote{There are actually two choices for this flag variety; we refer to \cite{CaraianiScholze} for a discussion of which one to choose.} Fix a prime $p$ and a prime $\mathfrak p$ of $E$ dividing $p$.

\begin{enumerate}
\item[{\rm (i)}] For any compact open subgroup $K^p\subset G(\A_f^p)$, there is a unique perfectoid space $\Sh_{K^p}$ over $E_{\mathfrak p}$ such that
\[
\Sh_{K^p}\sim \varprojlim_{K_p\subset G(\Q_p)} (\Sh_{K_pK^p}\otimes E_{\mathfrak p})^{\mathrm{rig}}\ .
\]
\item[{\rm (ii)}] There is a $G(\Q_p)$-equivariant Hodge-Tate period map
\[
\pi_{HT}: \Sh_{K^p}\to \Fl_{G,\mu}\ ,
\]
where we consider the right-hand side as an adic space over $E_{\mathfrak p}$.
\item[{\rm (iii)}] There is a Newton stratification
\[
\Fl_{G,\mu} = \bigsqcup_{b\in B(G,\mu)} \Fl_{G,\mu}^b
\]
into locally closed strata.
\item[{\rm (iv)}] If the Shimura variety is compact and of PEL type, then if $\overline{x}\in \Fl_{G,\mu}^b$ is a geometric point, the fibre $\pi_{HT}^{-1}(\overline{x})$ is the canonical lift of the perfection of the Igusa variety associated with $b$.
\end{enumerate}
\end{theorem}

We note that the Newton strata in (iii) are only defined on the adic space: sometimes they are nonempty but have no classical points!

\section{Applications to Langlands reciprocity}

The geometry of the Hodge-Tate period map has been used to obtain new results on the Langlands conjectures relating automorphic forms and Galois representations, especially in the case of torsion coefficients.

Let us first recall the results obtained in \cite{ScholzeTorsion}. For any reductive group $G$ over $\Q$ and a congruence subgroup $\Gamma\subset G(\Q)$, one can look at the locally symmetric $X_\Gamma= \Gamma\backslash X$, where $X$ is the symmetric space for $G(\R)$. To study Hecke operators, it is more convenient to switch to the adelic formalism, and consider
\[
X_K = G(\Q)\backslash (X\times G(\A_f)/K)
\]
for a compact open subgroup $K\subset G(\A_f)$, assumed sufficiently small from now on. The cohomology $H^i(X_K,\C)$ with complex coefficients can be computed in terms of automorphic forms. By the Langlands correspondence, one expects associated Galois representations. Much progress was made on these questions in case $X$ is a hermitian symmetric space, so that the $X_K$ are algebraic varieties over number fields, and their \'etale cohomology gives Galois representations. However, the groups $G=\GL_n$ for $n>2$ are not of this type.

It was conjectured by Grunewald in the 70's and later more precisely Ash, \cite{Ash}, that this relation between cohomology and Galois representations extends to the full integral cohomology groups $H^i(X_K,\Z)$, including their torsion subgroups, which can be enormous, especially in the case of hyperbolic $3$-manifolds, cf.~\cite{BergeronVenkatesh}.

\begin{theorem}[\cite{ScholzeTorsion}]\label{thm:torsion} Assume that $G=\Res_{F/\Q} \GL_n$ for some totally real or CM field $F$. Consider the abstract Hecke algebra $\mathbb T$ acting on $X_K$, generated by Hecke operators at good primes, and let $\mathbb T_K\subset \End(\bigoplus_i H^i(X_K,\Z))$ be the image of $\mathbb T$. For any maximal ideal $\mathfrak m\subset \mathbb T_K$, there is a continuous semisimple representation
\[
\overline{\rho}_{\mathfrak m}: \Gal(\overline{F}/F)\to \GL_n(\mathbb T_K/\mathfrak m)
\]
which is unramified at good primes, with Frobenius eigenvalues determined in terms of Hecke operators. If $\overline{\rho}_{\mathfrak m}$ is absolutely irreducible, then there is a nilpotent ideal $I_{\mathfrak m}\subset \mathbb T_{K,\mathfrak m}$ in the $\mathfrak m$-adic completion $\mathbb T_{K,\mathfrak m}$ of $\mathbb T_K$ whose nilpotence degree is bounded in terms of $n$ and $[F:\Q]$, such that there is a continuous representation
\[
\rho_{\mathfrak m}: \Gal(\overline{F}/F)\to \GL_n(\mathbb T_{K,\mathfrak m}/I_{\mathfrak m})
\]
which is unramified at good primes, with the characteristic polynomials of Frobenius elements determined in terms of Hecke operators.

In particular, for all cohomological automorphic representations $\pi$ of $G$, there exists a corresponding continuous semisimple Galois representation $\rho_\pi: \Gal(\overline{F}/F)\to \GL_n(\overline{\Q}_\ell)$ unramified at good primes.
\end{theorem}

The final part of this theorem was proved previously by Harris-Lan-Taylor-Thorne, \cite{HLTT}. The general strategy is to realize $X_K$ as a boundary component of the Borel-Serre compactification of a Shimura variety $\tilde{X}_{\tilde{K}}$ associated with a quasisplit symplectic or unitary group, and use the known existence of Galois representations for cusp forms on that space. The key part of the argument then is to show that all torsion cohomology classes on $\tilde{X}_{\tilde{K}}$ can be lifted to characteristic $0$, which is done using certain subtle results from $p$-adic Hodge theory in \cite{ScholzeTorsion}.

For applications, one needs a precise understanding of the behaviour of the Galois representations at bad primes as well, in particular for primes dividing $\ell$. To attack such questions, a better understanding of the torsion in $\tilde{X}_{\tilde{K}}$ is necessary. This has been obtained recently in joint work with Caraiani.

\begin{theorem}[\cite{CaraianiScholze2}]\label{thm:genericvanishing} Let $\widetilde{\mathbb T}$ be the Hecke algebra acting on the cohomology of $\tilde{X}_{\tilde{K}}$. Assume that the maximal ideal $\widetilde{\mathfrak m}\subset \widetilde{\mathbb T}$ is generic. Then
\[
H_c^i(\tilde{X}_{\tilde{K}},\Z_\ell)_{\widetilde{\mathfrak m}}=0
\]
for $i>d=\dim_\C \tilde{X}_{\tilde{K}}$.
\end{theorem}

Essentially by Poincar\'e duality, this implies that also $H^i(\tilde{X}_{\tilde{K}},\Z_\ell)_{\widetilde{\mathfrak m}}=0$ for $i<d$, and that $H^d(\tilde{X}_{\tilde{K}},\Z_\ell)_{\widetilde{\mathfrak m}}$ is $\ell$-torsion free. From the realization of $X_K$ in the boundary of $\tilde{X}_{\tilde{K}}$, one gets a long exact sequence
\[
\ldots\to H^d(\tilde{X}_{\tilde{K}},\Z_\ell)_{\widetilde{\mathfrak m}}\to H^d((X_K)^\prime,\Z_\ell)_{\widetilde{\mathfrak m}}\to H_c^{d+1}(\tilde{X}_{\tilde{K}},\Z_\ell)_{\widetilde{\mathfrak m}}=0\to \ldots\ ,
\]
where $(X_K)^\prime$ is a torus bundle over $X_K$, and the source is $\ell$-torsion free. Now the general strategy is to move any cohomology class in $H^i(X_K,\Z_\ell)_{\mathfrak m}$ to $H^d((X_K)^\prime,\Z_\ell)_{\widetilde{\mathfrak m}}$ by using the torus bundle to shift cohomological degrees, lift it to $H^d(\tilde{X}_{\tilde{K}},\Z_\ell)_{\widetilde{\mathfrak m}}$, and then use that this group injects into its rationalization, which can be expressed in terms of automorphic forms on $\tilde{G}$.

Let us comment on the proof of Theorem~\ref{thm:genericvanishing}. The result works for more general Shimura varieties,\footnote{The compact case has appeared in \cite{CaraianiScholze}.} so let us change the notation. Let $(\Sh_K)_{K\subset G(\A_f)}$ be a Shimura variety corresponding to a reductive group $G$ over $\Q$ together with some extra data, including a conjugacy class of minuscule cocharacters $\mu: \G_m\to G$, defined over the reflex field $E$. Fix a prime $p\neq \ell$ so that a certain genericity condition on the fixed maximal ideal of the Hecke algebra holds true at $p$, and a (sufficiently small) tame level $K^p\subset G(\A_f^p)$. We get (the minimal compactification of) the perfectoid Shimura variety at infinite level,
\[
\Sh_{K^p}^\ast\sim \varprojlim_{K_p} (\Sh_{K_p K^p}^\ast\otimes E_{\mathfrak p})^{\mathrm{rig}}\ ,
\]
where $\mathfrak p$ is a prime of $E$ dividing $p$. We have the Hodge-Tate period map
\[
\pi_{HT}: \Sh_{K^p}^\ast\to \Fl_{G,\mu}\ .
\]
The strategy now is to rewrite $R\Gamma_c(\Sh_{K^p},\Z_\ell)$ as $R\Gamma(\Fl_{G,\mu},R\pi_{HT\ast} j_! \Z_\ell)$, where $j: \Sh_{K^p}\to \Sh_{K^p}^\ast$ is the open immersion. The Hecke operators away from $p$ act trivially on the flag variety, so one can also rewrite
\[
R\Gamma_c(\Sh_{K^p},\Z_\ell)_{\mathfrak m} = R\Gamma(\Fl_{G,\mu},(R\pi_{HT\ast} j_! \Z_\ell)_{\mathfrak m})\ .
\]
The task is now to understand the sheaf $(R\pi_{HT\ast} j_! \Z_\ell)_{\mathfrak m}$. The first observation is that, with a suitable definition, it lies in $^p D^{\leq d}$ for the perverse $t$-structure; this uses that $R\pi_{HT\ast}$ is simultaneously affine and partially proper (but still has fibres of positive dimension -- a phenomenon only possible in this highly nonnoetherian setup). The other observation is that its fibres are given by the cohomology of Igusa varieties, by using Theorem~\ref{thm:fibreigusa}. The cohomology of Igusa varieties has been computed by Shin, \cite{Shin}, but only with $\overline{\Q}_\ell$-coefficients and in the Grothendieck group. Under the genericity condition, one finds that this always gives zero except if the point lies in $\Fl_{G,\mu}(\Q_p)\subset \Fl_{G,\mu}$. One can now play off these observations, which gives the conclusion that $(R\pi_{HT\ast} j_! \Z_\ell)_{\mathfrak m}$ is concentrated on the $0$-dimensional space $\Fl_{G,\mu}(\Q_p)$. Thus, there is no higher cohomology, and the bound $(R\pi_{HT\ast} j_! \Z_\ell)_{\mathfrak m}\in {}^p D^{\leq d}$ implies that $R\Gamma_c(\Sh_{K^p},\Z_\ell)_{\mathfrak m} = R\Gamma(\Fl_{G,\mu},(R\pi_{HT\ast} j_! \Z_\ell)_{\mathfrak m})$ is in degrees $\leq d$.

One nice aspect of this strategy is that it describes the cohomology of Shimura varieties in terms of the cohomology of certain sheaves on the flag varieties $\Fl_{G,\mu}$, and it becomes an interesting question to understand those sheaves themselves. This leads to a relation to the geometrization of the local Langlands correspondence conjectured by Fargues, cf.~Section~\ref{sec:twistor}.

Finally, let us mention that these results have led to the following applications.

\begin{theorem}[\cite{10Authors}]\label{thm:10authors} Let $F$ be a CM field.
\begin{enumerate}
\item For any elliptic curve $E$ over $F$, the $L$-function $L(E,s)$ has meromorphic continuation to $\C$. Moreover, $E$ satisfies the Sato-Tate conjecture.
\item For any cuspidal automorphic representation $\pi$ of $\GL_2(\A_F)$ whose archimedean component is of parallel weight $2$, the Ramanujan-Petersson conjecture holds true at all places.
\end{enumerate}
\end{theorem}

Theorem~\ref{thm:10authors} is proved by establishing the potential automorphy of $E$ and of all symmetric powers of $E$ resp.~$\pi$; here ``potential'' means after base change to (many) extensions $\tilde{F}/F$, as in \cite{BLGHT} where the case of totally real $F$ is proved. Thus, the second part follows Langlands's strategy for establishing the Ramanujan-Petersson conjecture in general, and it is the first time that this conjecture has been proved in a case where the associated motive (or some closely related motive) is not known to exist, and in particular the proof does not invoke Deligne's theorem on the Weil conjectures.

For the proof of Theorem~\ref{thm:10authors}, one wants to see that ``all'' Galois representations arise via Theorem~\ref{thm:torsion}. This follows Wiles' strategy \cite{Wiles}; more precisely, we use the variant proposed by Calegari-Geraghty, \cite{CalegariGeraghty}.

\section{$p$-adic twistor theory}\label{sec:twistor}

Further developments in $p$-adic geometry arose from the realization that the structures arising from the cohomology of proper smooth rigid-analytic varieties over the algebraically closed $p$-adic field $C$ can be naturally organized into a modification of vector bundles on the Fargues-Fontaine curve in a way closely resembling a reinterpretation of Hodge theory in terms of vector bundles on the twistor-$\mathbb P^1$.

Let us first recall the statements over $\C$. Consider the twistor-$\mathbb P^1$, i.e.~the nonsplit real form $\widetilde{\mathbb P^1_\R}$ of $\mathbb P^1$, which we will take to be given as the descent of $\mathbb P^1_\C$ to $\R$ via $z\mapsto -\frac 1{\overline{z}}$. We fix the point $\infty\in \widetilde{\mathbb P^1_\R}(\C)$ corresponding to $\{0,\infty\}\subset \mathbb P^1_\C$. There is an action of the nonsplit real torus $U(1)$ (that we consider as an algebraic group) on $\widetilde{\mathbb P^1_\R}$ fixing $\infty$.

\begin{proposition}[Simpson, {\cite[Section 5]{SimpsonHodgeFiltrNonabelian}}] The category of $U(1)$-equivariant semistable vector bundles on $\widetilde{\mathbb P^1_\R}$ is equivalent to the category of pure $\R$-Hodge structures.
\end{proposition}

Let us briefly recall the proof of this result. As the action of the algebraic group $U(1)$ on $\widetilde{\mathbb P^1_\R}\setminus \{\infty\}$ is simply transitive (as after base change to $\C$, it is the simply transitive action of $\G_m$ on itself), $U(1)$-equivariant vector bundles on $\widetilde{\mathbb P^1_\R}\setminus \{\infty\}$ are equivalent to $\R$-vector spaces $V$. Thus, given a $U(1)$-equivariant vector bundle $\mathcal E$ on $\widetilde{\mathbb P^1_\R}$, we get an $\R$-vector space $V$ such that
\[\mathcal E|_{\widetilde{\mathbb P^1_\R}\setminus \{\infty\}}\cong V\otimes_\R \OO_{\widetilde{\mathbb P^1_\R}\setminus \{\infty\}}
\]
equivariantly for the $U(1)$-action. Identifying the completion of $\widetilde{\mathbb P^1_\R}$ at $\infty$ with $\Spf \C[[t]]$, we get at the completion at $\infty$ a $U(1)$-equivariant $\C[[t]]$-lattice $\Lambda\subset V\otimes_\R \C((t))$. However, this is equivalent to a decreasing filtration $\Fil^i V_\C\subset V_\C$, via the Rees construction
\[
\Lambda = \sum_{i\in \Z} t^{-i}(\Fil^i V_\C)[[t]]\subset V\otimes_\R \C((t))\ .
\]
One then checks that $\mathcal E$ is semistable precisely when $(V,\Fil^i V_\C)$ defines a pure Hodge structure.

Thus, a description of the inverse functor is that one starts with a pure $\R$-Hodge structure $(V,\Fil^i V_\C)$ and the trivial vector bundle $\mathcal E^\prime = V\otimes_\R \OO_{\widetilde{\mathbb P^1_\R}}$, and modifies it at $\infty$ via the lattice $\Lambda$ defined above to obtain a new vector bundle $\mathcal E$ which is still $U(1)$-equivariant.

A very similar formalism exists in the $p$-adic case. Associated with the algebraically closed $p$-adic field $C$, there is the Fargues-Fontaine curve $\FF_C$, which is a regular noetherian scheme over $\Q_p$ of Krull dimension $1$ with a distinguished point $\infty\in \FF_C$. The completed local ring of $\FF_C$ at $\infty$ is Fontaine's period ring $B_\dR^+$, a complete discrete valuation ring with residue field $C$ and fraction field $B_\dR$.

\begin{theorem}[Fargues-Fontaine, \cite{FarguesFontaine}] All residue fields of $\FF_C$ at closed points are algebraically closed nonarchimedean extensions of $\Q_p$. Any vector bundle on $\FF_C$ is a direct sum of stable vector bundles, and there is a unique stable vector bundle $\OO_{\FF_C}(\lambda)$ for every rational slope $\lambda=\frac{s}{r}\in \Q$, which is of rank $r$ and degree $s$ (if $r$ and $s$ are chosen coprime with $r>0$).
\end{theorem}

Note that a similar result holds true for $\widetilde{\mathbb P^1_{\R}}$, but only half-integral slopes $\lambda\in \tfrac 12\Z$ occur in that case.

Given a proper smooth rigid-analytic space $X$ over $C$, one can form the trivial vector bundle
\[
\mathcal E = H^i_\et(X,\Z_p)\otimes_{\Z_p} \OO_{\FF_C}\ .
\]

\begin{theorem}[{\cite[Theorem 13.1, Theorem 13.8]{BMS}}]\label{thm:BdRpluscohom} There is a functorial $B_\dR^+$-lattice
\[
\Xi = H^i_\mathrm{crys}(X/B_\dR^+)\subset H^i_\et(X,\Z_p)\otimes_{\Z_p} B_\dR\ .
\]
If $X=X_0\hat{\otimes}_K C$ for some discretely valued subfield $K\subset C$, then $H^i_\mathrm{crys}(X/B_\dR^+) = H^i_\dR(X_0/K)\otimes_K B_\dR^+$ and the inclusion comes from the de~Rham comparison isomorphism.
\end{theorem}

Thus, one can form a modification $\mathcal E^\prime$ of $\mathcal E$ along $\infty$. The new vector bundle will in general be related to the (log-)crystalline cohomology of a (log-)smooth formal model. In other words, the comparison isomorphism between \'etale and crystalline cohomology can be understood in terms of a modification of vector bundles on the Fargues-Fontaine curve.

In fact, contrary to \'etale cohomology that is infinite-dimensional for spaces like the unit disc, it should be possible to define this modified vector bundle much more generally, without properness and without using formal models. The following conjecture arose from discussions of the author with Arthur-C{\'e}sar le Bras, and is related to results of Colmez--Nizio\l{}, \cite{ColmezNiziol}, \cite{Niziol}.

\begin{conjecture} There is a cohomology theory $H^i_{\FF_C}(X)$ for quasicompact separated smooth rigid spaces $X$ over $C$ taking values in vector bundles on $\FF_C$. If $X$ has an overconvergent model $X^\dagger$ (for example, if $X$ is affinoid), then the fibre of $H^i_{\FF_C}(X)$ at $\infty$ is the overconvergent de~Rham cohomology $H^i_{\mathrm{dR}}(X^\dagger/C)$ in the sense of Gro{\ss}e-Kl{\"o}nne, \cite{GrosseKloenne}, and the completion of $H^i_{\FF_C}(X)$ at $\infty$ is the overconvergent crystalline cohomology $H^i_{\mathrm{crys}}(X^\dagger/B_{\mathrm{dR}}^+)$ defined following \cite[Section 13]{BMS}.
\end{conjecture}

In other words, not only does de~Rham cohomology lift naturally to $B_\dR^+$ along the surjection $B_\dR^+\to C$ as can be explained in terms of a version of crystalline cohomology, but it does actually deform into a vector bundle on all of $\FF_C$.

Now we define for any local Shimura datum $(G,b,\mu)$ the local Shimura variety $\mathcal M_{(G,b,\mu),K}$ in terms of modifications of $G$-torsors on the Fargues-Fontaine curve. 

\begin{theorem}[Fargues, \cite{FarguesGTorsor}]\label{thm:BGCurve} For any reductive group $G$ over $\Q_p$, there is a natural bijection $b\mapsto \mathcal E_b$ between Kottwitz' set $B(G)$ of $\sigma$-conjugacy classes and the set of $G$-torsors on $\FF_C$ up to isomorphism. Moreover, $b\in B(G,\mu)\subset B(G)$ if and only if $\mathcal E_b$ can be written as a modification at $\infty$ of type $\mu$ of the trivial $G$-torsor $\mathcal E_1$.
\end{theorem}

The final statement appears in \cite{RapoportAppendix}. Now we define
\[
\mathcal M_{(G,b,\mu),\infty}(C) = \varprojlim_K \mathcal M_{(G,b,\mu),K}(C)
\]
to be the set of all modifications $\mathcal E_1\dashrightarrow \mathcal E_b$ at $\infty$ of type $\mu$. The group $G(\Q_p)$ is the group of automorphisms of $\mathcal E_1$, and thus acts on this inverse limit; then for all $K\subset G(\Q_p)$,
\[
\mathcal M_{(G,b,\mu),K}(C) = \mathcal M_{(G,b,\mu),\infty}(C)/K\ .
\]
In order to make $\mathcal M_{(G,b,\mu),K}$ into a rigid space, we need to define it as a moduli problem. For the rest of this section, we make extensive use of the theory of perfectoid spaces, and assume that the reader is familiar with it, cf.~\cite{ScholzeICM}; some of the structures that appear now will however be motivated in the next section. We use that for any perfectoid space $S$ of characteristic $p$ together with an untilt $S^\sharp$ over $\Q_p$, one can construct a relative Fargues-Fontaine curve $\FF_S$ that is an adic space over $\Q_p$ and comes with a section $\infty: \Spa S^\sharp\to \FF_S$. Then $\mathcal M_{(G,b,\mu),\infty}(S)$ parametrizes modifications $\mathcal E_1|_{\FF_S}\dashrightarrow \mathcal E_b|_{\FF_S}$ at $\infty$ of type $\mu$ as before. This defines a structure of a diamond, \cite{ScholzeDiamonds}:

\begin{definition} Let $\Perf$ be the category of perfectoid spaces of characteristic $p$. A diamond is a pro-\'etale sheaf $Y$ on $\Perf$ that can be written as a quotient $Y=X/R$ of a perfectoid space $X$ by a pro-\'etale equivalence relation $R\subset X\times X$.
\end{definition}

An example of a diamond is given by the sheaf $\Spd \Q_p$ that attaches to any perfectoid space $S$ of characteristic $p$ the set of all untilts $S^\sharp$ over $\Q_p$. More generally, if $X$ is an adic space over $\Q_p$, one can define a diamond $X^\diamond$ whose $S$-valued points are given by an untilt $S^\sharp$ over $\Q_p$ together with a map $S^\sharp\to X$.

\begin{theorem}[\cite{Berkeley}] For any nonarchimedean field $L/\Q_p$, the functor $X\mapsto X^\diamond$ defines a fully faithful functor from the category of seminormal rigid spaces over $L$ to the category of diamonds over $\Spd L=(\Spa L)^\diamond$.

The diamond $\mathcal M_{(G,b,\mu),K} = \mathcal M_{(G,b,\mu),\infty} / \underline{K}$ over $\Spd \breve{E}$ is the image of a smooth rigid space over $\breve{E}$ under this functor.
\end{theorem}

The proof makes use of the results of Kedlaya, \cite{KedlayaICM}, and Kedlaya-Liu, \cite{KedlayaLiu}, on families of vector bundles on the Fargues-Fontaine curve.

From here, it becomes natural to consider much more general spaces, parametrizing modifications of arbitrary $G$-bundles at several points with modifications of arbitrary type, not necessarily minuscule. Such spaces live over the base space $\Spd \Q_p\times \ldots \times \Spd \Q_p$ parametrizing the points of modification: These products give good meaning to the non-existent products $\Spec \Q_p\times_{\Spec \F_1}\ldots\times_{\Spec \F_1} \Spec \Q_p$. The resulting spaces are now in general just diamonds. The purpose of the foundational manuscript \cite{ScholzeDiamonds} is to develop a solid theory of \'etale cohomology for diamonds with the aim of using these spaces to obtain a general local Langlands correspondence.

These spaces are naturally organized into Hecke stacks acting on the stack $\mathrm{Bun}_G$ of $G$-bundles on the Fargues-Fontaine curve. Fargues realized that this gives rise to a picture perfectly resembling the geometric Langlands correspondence for classical smooth projective curves, \cite{FarguesSurvey}, \cite{FarguesICM}, \cite{FarguesScholze}. This in particular involves perverse $\ell$-adic sheaves on $\mathrm{Bun}_G$ whose pullback to the flag varieties $\Fl_{G,\mu}$ should arise globally from the construction $R\pi_{HT\ast} \Z_\ell$ sketched in the previous section.

It is too early to state the expected theorems about the applications to the local Langlands correspondence, but let us mention that recently Kaletha-Weinstein, \cite{KalethaWeinstein}, have used the \'etale cohomology of diamonds as developed in \cite{ScholzeDiamonds} to prove Kottwitz' conjecture, \cite{RapoportICM}, about the realization of the Jacquet-Langlands correspondence in the cohomology of local Shimura varieties.

\begin{remark} As the theory of diamonds is critical to the foundations, as is the possibility of defining $\Spd \Q_p\times \Spd \Q_p$, let us quickly give a description of the product $\Spd \Q_p\times \Spd \Q_p$. Consider the open unit disc $\mathbb D_{\Q_p} = \{z\mid |z|<1\}$ as embedded in the multiplicative group via $z\mapsto 1+z$; then $\mathbb D_{\Q_p}$ itself is a group object. Let $\widetilde{\mathbb D}_{\Q_p}$ be the inverse limit of $\mathbb D_{\Q_p}$ along $z\mapsto (1+z)^p-1$. Then $\widetilde{\mathbb D}_{\Q_p}$ is (pre)perfectoid, and admits a natural $\Q_p$-action. Now
\[
\Spd \Q_p\times \Spd \Q_p = (\widetilde{\mathbb D}_{\Q_p}\setminus \{0\})^\diamond/\underline{\Z_p^\times}\ .
\]
On the right-hand side, one factor of $\Spd \Q_p$ arises from the structure map $\widetilde{\mathbb D}_{\Q_p}\to \Spa \Q_p$, but the other factor gets realized in terms of the perfectoid punctured unit disc, and thus has become geometric.

These constructions in fact lead to a description of the absolute Galois group of $\Q_p$ as a geometric fundamental group:

\begin{theorem}[\cite{WeinsteinGeomGal}] For any algebraically closed nonarchimedean field $C/\Q_p$, the \'etale fundamental group of $(\widetilde{\mathbb D}_C\setminus \{0\})/\underline{\Q_p^\times}$ agrees with the absolute Galois group of $\Q_p$.
\end{theorem}

A search for a hypothetical space $\Spec \Z\times \Spec \Z$ thus seems closely related to a realization of the absolute Galois group of $\Q$ as a geometric fundamental group. For a step in this direction, let us mention the following result, which uses the ring of rational Witt vectors $W_{\mathrm{rat}}(R)\subset W_{\mathrm{big}}(R)$. For all $r\in R$, there is the Teichm\"uller lift $[r]\in W_{\mathrm{rat}}(R)$.

\begin{theorem}[\cite{KucharczykScholze}]\label{thm:Wrat} Let $L$ be a field of characteristic $0$ that contains all roots of unity, and fix an embedding $\Q/\Z\hookrightarrow L^\times$, $1/n\mapsto \zeta_n$. Then the category of finite extensions of $L$ is equivalent to the category of connected finite coverings of the topological space
\[
X(L)\subset (\Spec W_{\mathrm{rat}}(L))(\C)
\]
that is the connected component singled out by the condition that $[\zeta_n]$ maps to $e^{2\pi i/n}\in \C$ for all $n\geq 1$.
\end{theorem}

In fact, $X(L)$ has a deformation retract to a compact Hausdorff space, which gives a realization of the absolute Galois group of $L$ as the profinite fundamental group of a compact Hausdorff space! This formally implies that the absolute Galois group of $L$ is torsion-free, cf.~\cite[Proposition 7.10]{KucharczykScholze}.

The space $X(L)$ gives rise to certain non-profinitely complete structures on natural arithmetic invariants. If $L$ is the cyclotomic extension of $\Q$, one can show that the actual fundamental group of $X(L)$ given by topological loops is a proper dense subgroup of the absolute Galois group. Moreover, it acts naturally on the set $\log \overline{\Q}\subset \C$ of logarithms of algebraic numbers, compatibly with $\exp: \log \overline{\Q}\to \overline{\Q}^\times$. Similarly, for general $L$ as in the theorem, the \v{C}ech cohomology groups $H^i(X(L),\Z)$ give certain torsion-free non-profinitely completed abelian groups with
\[
H^i(X(L),\Z)/n = H^i(\Gal(\overline{L}/L),\Z/n\Z)
\]
for all $n\geq 1$, cf.~\cite[Theorem 1.8]{KucharczykScholze}.
\end{remark}

\section{$p$-adic shtukas}\label{sec:padicshtuka}

Going further, one can reinterpret the linear-algebraic structures that appeared in the last section in terms of a $p$-adic analogue of shtukas. Let us first recall the basic definition of a shtuka from the function field case. There is a version of the following definitions for any reductive group $G$ obtained by replacing vector bundles by $G$-torsors.

\begin{definition} Let $X$ be a smooth projective curve over $\F_p$ and let $S$ be a scheme over $\F_p$. A shtuka over $S$ relative to $X$ with legs at $x_1,\ldots,x_n: S\to X$ is a vector bundle $\mathcal E$ over $S\times_{\F_p} X$ together with an isomorphism 
\[
(\Frob_S^\ast \mathcal E)|_{S\times_{\F_p} X\setminus \bigcup_{i=1}^n \Gamma_{x_i}}\cong \mathcal E|_{S\times_{\F_p} X\setminus \bigcup_{i=1}^n \Gamma_{x_i}}\ ,
\]
where $\Gamma_{x_i}: S\hookrightarrow S\times_{\F_p} X$ is the graph of $x_i$.
\end{definition}

Drinfeld used moduli spaces of shtukas with two legs to prove the global Langlands correspondence for $\GL_2$, \cite{DrinfeldGL2}. This was generalized to $\GL_n$, still using moduli spaces with two legs, by L.~Lafforgue, \cite{LLafforgue}. Recently, V.~Lafforgue has used all moduli spaces with an arbitrary number of legs simultaneously to prove the automorphic to Galois direction of the global Langlands correspondence for any reductive group $G$, \cite{VLafforgue}.

There is a corresponding notion of local shtuka, where one works with the local curve $X=\Spf \F_p[[t]]$. The legs are now parametrized by maps $x_1,\ldots,x_n: S\to X=\Spf \F_p[[t]]$, i.e.~locally nilpotent elements $t_i\in \OO_S(S)$. As this contains no topological information, we pass to the world of rigid geometry. For simplicity, we discuss the case of geometric points, and so assume that $S=\Spa C$ for some complete algebraically closed nonarchimedean field $C$ of characteristic $p$; we will later allow more general base spaces again. The legs are maps $S\to X$, which are parametrized by topologically nilpotent elements $t_i\in C$, i.e.~$|t_i|<1$. The fibre product $S\times_{\F_p} X$ becomes the open unit ball $\mathbb D_C=\{t\in C\mid |t|<1\}$, and the legs give rise to points $x_i\in \mathbb D_C$. There is a Frobenius $\Frob_C$ acting on $\mathbb D_C$, coming from the Frobenius on $C$. This is not a map of rigid spaces over $C$, but it does exist in the category of adic spaces.

\begin{definition} A (local) shtuka over $S=\Spa C$ relative to $X=\Spf \F_p[[t]]$ with legs at $x_1,\ldots,x_n: S\to X$ given by elements $t_i\in C$, $|t_i|<1$, is a vector bundle $\mathcal E$ over $\mathbb D_C = S\times_{\F_p} X$ together with an isomorphism
\[
(\Frob_C^\ast \mathcal E)|_{\mathbb D_C\setminus \{t_1,\ldots,t_n\}}\cong \mathcal E|_{\mathbb D_C\setminus \{t_1,\ldots,t_n\}}
\]
that is meromorphic along the $t_i$.
\end{definition}

The main observation of \cite{Berkeley} is that it is possible to give a mixed-charac\-te\-ris\-tic version of this definition, where $X=\Spf \Z_p$. As a geometric point, one still takes $S=\Spa C$ where $C$ is a complete algebraically closed nonarchimedean field of characteristic $p$. Naively, there are now no interesting maps $S=\Spa C\to X=\Spf \Z_p$; there is exactly one, which factors over $\Spec \F_p$. We will see momentarily how to solve this problem. But let us first face the other problem of defining a reasonable fibre product $``S\times X"$, where the product should be over $\Spec \F_p$. The basic insight is that for any perfect ring $R$, $``\Spec R\times \Spf \Z_p"$ should be given by $\Spf W(R)$ with the $p$-adic topology on $W(R)$. If $R$ is itself an adic ring with the $I$-adic topology, then $``\Spf R\times \Spf \Z_p = \Spf W(R)"$, where $W(R)$ is equipped with the $(p,[I])$-adic topology.

In particular, we should set $``\Spf \OO_C\times \Spf \Z_p" = \Spf(W(\OO_C))$. The ring $W(\OO_C)$ plays an important role in $p$-adic Hodge theory, and is traditionally called $A_\inf$ there. If $\varpi\in C$ is a topologically nilpotent unit, then passing from $\Spf \OO_C$ to its generic fibre $\Spa C$ amounts to inverting $\varpi$. Thus, we define the open subspace
\[
``\Spa C\times \Spa \Z_p" = \{[\varpi]\neq 0\}\subset \Spa(W(\OO_C))\ ,
\]
which is a well-defined adic space. The following proposition shows that this behaves very much like a classical smooth rigid curve, except that it does not live over any base field.

\begin{proposition}[\cite{FarguesFontaine}, \cite{KedlayaNoetherian}] For any connected open affinoid subspace
\[
U=\Spa(R,R^+)\subset ``\Spa C\times \Spa \Z_p"\ ,
\]
the ring $R$ is a principal ideal domain. For any maximal ideal $\mathfrak m\subset R$, the quotient $C_{\mathfrak m} = R/\mathfrak m$ is a complete algebraically closed nonarchimedean field. Moreover, there is a canonical isomorphism
\[
\varprojlim_{x\mapsto x^p} \OO_{C_{\mathfrak m}}/p\cong \OO_C\ ,
\]
i.e.~$C\cong C_{\mathfrak m}^\flat$ is the tilt of $C_{\mathfrak m}$ in the terminology of \cite{ScholzeThesis}. Conversely, for any complete algebraically closed nonarchimedean field $C^\sharp$ with an isomorphism $(C^\sharp)^\flat\cong C$, there is a unique maximal ideal $\mathfrak m\subset R$ (for any large enough open subset $U$) such that $C^\sharp=R/\mathfrak m$, compatibly with the identification of $(C^\sharp)^\flat$ with $C$.
\end{proposition}

In other words, ``classical points" of $``\Spa C\times \Spa \Z_p"$ parametrize untilts of $C$. This shows that we should reinterpret the data of the maps $x_i: S=\Spa C\to X=\Spa \Z_p$ as the data of untilts $C_1^\sharp,\ldots,C_n^\sharp$ of $C$.

\begin{definition} A (local/$p$-adic) shtuka over $S=\Spa C$ relative to $X=\Spf \Z_p$ with legs at $x_1,\ldots,x_n$ given by untilts $C_1^\sharp,\ldots,C_n^\sharp$ of $C$, is a vector bundle $\mathcal E$ over $``\Spa C\times \Spa \Z_p"$ together with an isomorphism
\[
(\Frob_C^\ast \mathcal E)|_{``\Spa C\times \Spa \Z_p"\setminus \{x_1,\ldots,x_n\}}\cong \mathcal E|_{``\Spa C\times \Spa \Z_p"\setminus \{x_1,\ldots,x_n\}}
\]
that is meromorphic along the $x_i$.
\end{definition}

In the case of one leg, these structures are closely related to the structures that appeared in the last section, by the following theorem of Fargues.

\begin{theorem}[Fargues, \cite{Berkeley}]\label{thm:fargues} Assume that $C^\sharp$ is an untilt over $\Q_p$, let $\infty: \Spec C^\sharp\to \FF_C$ be the corresponding point of the Fargues-Fontaine curve,\footnote{The Fargues-Fontaine curve depends only on the tilt of $C^\sharp$, but the point $\infty$ depends on $C^\sharp$.} and let $B_\dR^+ = B_\dR^+(C^\sharp)$ be the complete local ring at $\infty$ with quotient field $B_\dR$. The following categories are equivalent.
\begin{enumerate}
\item[{\rm (i)}] Shtukas over $S=\Spa C$ relative to $X=\Spf \Z_p$ with one leg at $C^\sharp$.
\item[{\rm (ii)}] Pairs $(T,\Xi)$, where $T$ is a finite free $\Z_p$-module and $\Xi\subset T\otimes_{\Z_p} B_\dR$ is a $B_\dR^+$-lattice.
\item[{\rm (iii)}] Quadruples $(\mathcal F_1,\mathcal F_2,\beta,T)$, where $\mathcal F_1$ and $\mathcal F_2$ are two vector bundles on the Fargues-Fontaine curve $\FF_C$, $\beta: \mathcal F_1|_{\FF_C\setminus \{\infty\}}\cong \mathcal F_2|_{\FF_C\setminus \{\infty\}}$ is an isomorphism, and $T$ is a finite free $\Z_p$-module such that $\mathcal F_1 = T\otimes_{\Z_p} \OO_{\FF_C}$ is the corresponding trivial vector bundle.
\item[{\rm (iv)}] Breuil-Kisin-Fargues modules, i.e.~finite free $A_\inf=W(\OO_C)$-modules $M$ together with a $\varphi$-linear isomorphism $M[\tfrac 1{\varphi^{-1}(\xi)}]\cong M[\tfrac 1{\xi}]$, where $\xi\in A_\inf$ is a generator of $\ker(A_\inf\to \OO_{C^\sharp})$.
\end{enumerate}
\end{theorem}

Note that local Shimura varieties are parametrizing data of type (iii). The equivalence with (i) shows that one can regard local Shimura varieties as moduli spaces of shtukas, and it gives a possibility of formulating an integral model for the local Shimura varieties that we will now discuss.

In \cite{Berkeley}, we construct moduli spaces of shtukas with any number of legs. However, already the case of one leg has important applications, so let us focus on this case. In addition to the group $G$ over $\Q_p$, we need to fix a model $\mathcal G$ of $G$ over $\Z_p$, and the most important case is when $\mathcal G$ is parahoric. To define the moduli problem, we need to bound the modification at the leg; this leads to a version of the affine Grassmannian. Here, $\Spd \Z_p$ denotes the functor on $\Perf$ taking any perfectoid space $S$ of characteristic $p$ to the set of untilts $S^\sharp$ over $\Z_p$. If $S=\Spa(R,R^+)$ is affinoid, then $S^\sharp=\Spa(R^\sharp,R^{\sharp+})$, and there is a natural ring $B_\dR^+(R^\sharp)$ that surjects onto $R^\sharp$ with kernel generated by some nonzerodivisor $\xi\in B_\dR^+(R^\sharp)$ such that $B_\dR^+(R^\sharp)$ is $\xi$-adically complete. This interpolates between the following cases:
\begin{enumerate}
\item If $R^\sharp=C^\sharp$ is an untilt over $\Q_p$ of an algebraically closed nonarchimedean field $R=C$ of characteristic $p$, then $B_\dR^+(C^\sharp)$ is Fontaine's ring considered previously.
\item If $R^\sharp=R$ is of characteristic $p$, then $B_\dR^+(R^\sharp)=W(R)$ is the ring of $p$-typical Witt vectors.
\end{enumerate}

\begin{definition} Let $\mathcal G$ be a parahoric group scheme over $\Z_p$. The Beilinson-Drin\-feld Grassmannian
\[
\Gr_{\mathcal G,\Spd \Z_p}\to \Spd \Z_p
\]
is the moduli problem on $\Perf$ taking an affinoid perfectoid space $S=\Spa(R,R^+)$ of characteristic $p$ to the set of untilts $S^\sharp=\Spa(R^\sharp,R^{\sharp+})$ over $\Z_p$ together with a $\mathcal G$-torsor over $B_\dR^+(R^\sharp)$ that is trivialized over $B_\dR(R^\sharp)$.
\end{definition}

In the following, we make use of the following proposition to identify perfect schemes or formal schemes as certain pro-\'etale sheaves. In both cases, one can define a functor $X\mapsto X^\diamond$, where $X^\diamond$ parametrizes untilts over $X$.

\begin{proposition}\label{prop:fullyfaithful} The functor $X\mapsto X^\diamond$ defines a fully faithful functor from the category of perfect schemes to the category of pro-\'etale sheaves on $\Perf$. Similarly, for any nonarchimedean field $L$ over $\Q_p$, the functor $X\mapsto X^\diamond$ from normal and flat formal schemes locally formally of finite type over $\Spf \OO_L$ to pro-\'etale sheaves over $\Spd \OO_L$ is fully faithful.
\end{proposition}

The proof of the second part relies on a result of Louren\c{c}o, \cite{Lourenco}, that recovers formal schemes as in the proposition from their generic fibre and the perfection of their special fibre, together with the specialization map.

\begin{theorem}[\cite{ZhuWitt}, \cite{BhattScholzeWitt}, \cite{Berkeley}] The special fibre of $\Gr_{\mathcal G,\Spd \Z_p}$ is given by the Witt vector affine Grassmannian $\Gr_{\mathcal G}^W$ constructed in \cite{ZhuWitt} and \cite{BhattScholzeWitt} that can be written as an increasing union of perfections of projective varieties along closed immersions. The generic fibre of $\Gr_{\mathcal G,\Spd \Z_p}$ is the $B_\dR^+$-affine Grassmannian of $G$, and can be written as an increasing union of proper diamonds along closed immersions.

For a conjugacy class of minuscule cocharacters $\mu$ defined over $E$, there is a natural closed immersion
\[
\Fl_{G,\mu}^\diamond\hookrightarrow \Gr_{\mathcal G,\Spd \Z_p}\times_{\Spd \Z_p} \Spd E\ .
\]
\end{theorem}

We let
\[
\mathbb M_{(\mathcal G,\mu)}^{\mathrm{loc}}\subset \Gr_{\mathcal G,\Spd \Z_p}\times_{\Spd \Z_p} \Spd \OO_E
\]
be the closure of $\Fl_{G,\mu}^\diamond$. We conjecture that it comes from a normal and flat projective scheme over $\Spec \OO_E$ under Proposition~\ref{prop:fullyfaithful}. This would give a group-theoretic definition of the local model in the theory of Shimura varieties, cf.~e.g.~\cite{PappasICM}.

\begin{definition} Given local Shimura data $(G,b,\mu)$ and a parahoric model $\mathcal G$, let
\[
\mathcal M^{\mathrm{int}}_{(\mathcal G,b,\mu)}\to \Spd \breve{\OO}_E
\]
be the moduli problem taking $S\in \Perf$ to the set of untilts $S^\sharp$ over $\breve{\OO}_E$ together with a $\mathcal G$-shtuka over $``S\times \Spa\Z_p"$ with one leg at $S^\sharp$ that is bounded by $\mathbb M_{(\mathcal G,\mu)}^{\mathrm{loc}}$, and a trivialization of the shtuka at the boundary of $``S\times \Spa \Z_p"$ by $\mathcal E_b$.
\end{definition}

For a precise formulation of the following result, we refer to \cite{Berkeley}.

\begin{theorem}[\cite{Berkeley}]\label{thm:rapoportzink} This definition recovers Rapoport-Zink spaces.
\end{theorem}

Using this group-theoretic characterization of Rapoport-Zink spaces, we can obtain new isomorphisms between different Rapoport-Zink spaces.

\begin{corollary} The conjectures of Rapoport-Zink, \cite{RapoportZinkAlternative}, and Kudla-Rapoport-Zink, \cite{KudlaRapoportZink}, on alternative descriptions of the Drinfeld moduli problem hold true.

In particular, by \cite{KudlaRapoportZink}, one gets an integral version and moduli-theoretic proof of \v{C}erednik's $p$-adic uniformization theorem, \cite{Cerednik}.
\end{corollary}

To prove Theorem~\ref{thm:rapoportzink}, one needs to see that $p$-adic shtukas are related to the cohomology of algebraic varieties or $p$-divisible groups. This is the subject of integral $p$-adic Hodge theory that we will discuss next.

\section{Integral $p$-adic Hodge theory}

The following question arises naturally from Theorem~\ref{thm:fargues}. As it is more natural in this section, let us change notation, and start with an algebraically closed $p$-adic field $C$ and consider shtukas over its tilt $C^\flat$ with one leg at the untilt $C$ of $C^\flat$. Given a proper smooth rigid-analytic space $X$ over $C$, consider $T=H^i_\et(X,\Z_p)/(\mathrm{torsion})$ with the finite free $B_\dR^+$-module $\Xi=H^i_{\mathrm{crys}}(X/B_\dR^+)\subset T\otimes_{\Z_p} B_\dR$ from Theorem~\ref{thm:BdRpluscohom}. By Theorem~\ref{thm:fargues}, there is a corresponding Breuil-Kisin-Fargues module $H^i_{A_\inf}(X)$; we normalize it here so that the Frobenius is an isomorphism after inverting $\xi$~resp.~$\varphi(\xi)$. Can one give a direct cohomological construction of this?

We expect that without further input, the answer is no; in fact, Fargues's equivalence is not exact, and does manifestly not pass to the derived category. However, in joint work with Bhatt and Morrow, \cite{BMS}, we show that it is possible once a proper smooth formal model $\mathfrak X$ of $X$ is given.\footnote{All of these results were recently extended to the case of semistable reduction by \v{C}esnavi\v{c}ius-Koshikawa, \cite{CesnaviciusKoshikawa}.}

\begin{theorem}[\cite{BMS}]\label{thm:bms} Let $\mathfrak X$ be a proper smooth formal scheme over $\Spf \OO_C$ with generic fibre $X$. There is a perfect complex $R\Gamma_{A_\inf}(\mathfrak X)$ of $A_\inf$-modules together with a $\varphi$-linear map $\varphi: R\Gamma_{A_\inf}(\mathfrak X)\to R\Gamma_{A_\inf}(\mathfrak X)$ that becomes an isomorphism after inverting $\xi$ resp.~$\varphi(\xi)$. Each $H^i_{A_\inf}(\mathfrak X)$ is a finitely presented $A_\inf$-module that becomes free over $A_\inf[\tfrac 1p]$ after inverting $p$. Moreover, one has the following comparison results.
\begin{enumerate}
\item[{\rm (i)}] Crystalline comparison: $R\Gamma_{A_\inf}(\mathfrak X)\buildrel{\mathbb L}\over\otimes_{A_\inf} W(k)\cong R\Gamma_{\mathrm{crys}}(\mathfrak X_k/W(k))$, $\varphi$-equivariantly, where $k$ is the residue field of $\OO_C$.
\item[{\rm (ii)}] De~Rham comparison: $R\Gamma_{A_\inf}(\mathfrak X)\buildrel{\mathbb L}\over\otimes_{A_\inf} \OO_C\cong R\Gamma_\dR(\mathfrak X/\OO_C)$.
\item[{\rm (iii)}] \'Etale comparison: $R\Gamma_{A_\inf}(\mathfrak X)\otimes_{A_\inf} W(C^\flat)\cong R\Gamma_\et(X,\Z_p)\otimes_{\Z_p} W(C^\flat)$, $\varphi$-equivariantly.
\end{enumerate}
Moreover, if $H^i_{\mathrm{crys}}(\mathfrak X_k/W(k))$ is $p$-torsion free\footnote{Equivalently, $H^i_\dR(\mathfrak X/\OO_C)$ is $p$-torsion free, cf.~\cite[Remark 14.4]{BMS}.}, then $H^i_{A_\inf}(\mathfrak X)$ is finite free over $A_\inf$, and agrees with $H^i_{A_\inf}(X)$ as defined above.\footnote{In that situation, part (iii) implies that also $H^i_\et(X,\Z_p)$ is $p$-torsion free.}
\end{theorem}

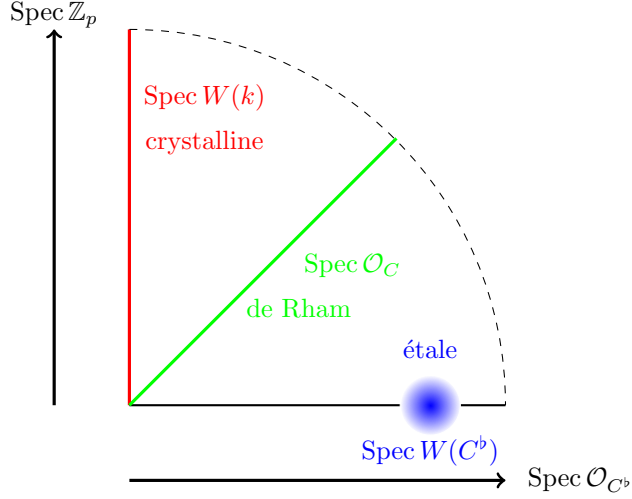
\begin{figure}
\centering

\begin{tikzpicture}[auto]
\draw[->, very thick] (-1,0) -- (-1,5);
\draw[->, very thick] (0,-1) -- (5,-1);
\draw (6,-1) node {$\Spec \OO_{C^\flat}$};
\draw (-1,5.2) node {$\Spec \Z_p$};
\draw[very thick, color=red] (0,0) -- (0,5);
\draw[thick] (0,0) -- (5,0);
\draw[very thick, color=green] (0,0) -- (3.55,3.55);
\draw[color=red] (1,3.5) node {$\mathrm{crystalline}$};
\draw[color=red] (1,4.1) node {$\Spec W(k)$};
\draw[dashed] (0,5) arc [start angle=90, end angle=0, radius=5];
\draw[inner color=blue, color=white] (4,0) circle [radius=.4];
\draw[color=blue] (4,-.6) node {$\Spec W(C^\flat)$};
\draw[color=blue] (4,0.75) node {$\mathrm{\acute{e}tale}$};
\draw[color=green] (32:3.5) node {$\Spec \OO_C$};
\draw[color=green] (30:2.6) node {$\mathrm{de~Rham}$};
\end{tikzpicture}

\caption{A picture of some parts of $\Spec A_{\inf}=``\Spec \OO_{C^\flat}\times \Spec \Z_p"$.}\label{fig:Ainf}
\end{figure}

The theorem implies a similar result for $p$-divisible groups, which is the key input into the proof of Theorem~\ref{thm:rapoportzink}. On the other hand, the theorem has direct consequences for the behaviour of torsion under specialization from characteristic $0$ to charateristic $p$.

\begin{corollary} For any $i\geq 0$, one has $\dim_k H^i_\dR(\mathfrak X_k/k)\geq \dim_{\F_p} H^i_\et(X,\F_p)$, and
\[
\mathrm{length}_{W(k)} H^i_{\mathrm{crys}}(\mathfrak X_k/W(k))_{\mathrm{tor}}\geq \mathrm{length}_{\Z_p} H^i_\et(X,\Z_p)_{\mathrm{tor}}\ .
\]
\end{corollary}

For example, for an Enriques surface $X$ over $C$, so that it has a double cover by a K3 surface, there is $2$-torsion in the second \'etale cohomology group. By the theorem, this implies that for all Enriques surfaces in characteristic $2$ (all of which deform to characteristic $0$), one has $2$-torsion in the second crystalline cohomology, or equivalently $H^1_\dR\neq 0$, contrary to the situation in any other characteristic. This is a well-known ``pathology'' of Enriques surfaces in characteristic $2$ that finds a natural explanation here. In this case, one actually has equality in the corollary; however, in \cite{BMS} we give an example of a projective smooth surface over $\Z_2$ whose \'etale cohomology is torsionfree while there is torsion in crystalline cohomology.

The construction of the $A_\inf$-cohomology given in \cite{BMS} was the end result of a long detour that led the author to study topological Hochschild and cyclic homology. This started with the paper \cite{HesselholtOC} of Lars Hesselholt that computes $\THH(\OO_C)$. This is an $E_\infty$-ring spectrum with an action by the circle group $\mathbb T = S^1$. In particular, one can form the homotopy fixed points $\TC^-(\OO_C) = \THH(\OO_C)^{h\mathbb T}$ to get another $E_\infty$-ring spectrum.

\begin{theorem}[Hesselholt, \cite{HesselholtOC}] The homotopy groups of the $p$-completion of $\TC^-(\OO_C)=\THH(\OO_C)^{h\mathbb T}$ are given by $A_\inf$ in all even degrees, and $0$ in all odd degrees. The generators in degree $2$ and $-2$ multiply to $\xi\in A_\inf$.
\end{theorem}

If now $X$ is a proper smooth (formal) scheme over $\OO_C$, it follows that the $p$-completion of $\TC^-(X)$ is a perfect module over the $p$-completion of $\TC^-(\OO_C)$, and in particular its homotopy groups are $A_\inf$-modules. Moreover, topological Hochschild homology comes with a $\mathbb T$-equivariant Frobenius operator $\THH(X)\to \THH(X)^{tC_p}$, where $C_p\subset \mathbb T$ is the cyclic group of order $p$ and $S^{tC_p} = \mathrm{cone}(\mathrm{Nm}: S_{hC_p}\to S^{hC_p})$ denotes the Tate construction. This gives rise to the desired Frobenius on $\TC^-$ by passing to homotopy $\mathbb T$-fixed points. In the classical formulation, $\THH$ has even more structure as encoded in the structure of a cyclotomic spectrum. However, in joint work with Nikolaus, we proved that this extra structure is actually redundant.

\begin{theorem}[\cite{NikolausScholze}]\label{thm:niksch} On respective subcategories of bounded below objects, the $\infty$-category of cyclotomic spectra (in the sense of Hesselholt-Madsen, \cite{HesselholtMadsen}, and made more precise by Blum\-berg-Mandell, \cite{BlumbergMandell}) is equivalent to the $\infty$-category of naively $\mathbb T$-equivariant spectra $S$ together with a $\mathbb T$-equivariant map $\varphi_p: S\to S^{tC_p}$ for all primes $p$.\footnote{No compatibility between the different $\varphi_p$ is required.}
\end{theorem}

The only issue with this construction of the $A_\inf$-cohomology theory (besides its heavy formalism) is that one gets an essentially $2$-periodic cohomology theory. To get the desired cohomology theory itself, one needs to find a filtration on $\TC^-(X)$ and then pass to graded pieces; this is similar to the ``motivic'' filtration on algebraic $K$-theory with graded pieces given by motivic cohomology, and the existence of such a filtration was first conjectured by Hesselholt. Using perfectoid techniques, it is possible to construct the desired filtration on $\TC^-$; this is the content of \cite{BMS2}. On the other hand, in the original paper \cite{BMS}, we were able to build the theory independently of $\THH$. We refer to \cite{BMS} and the surveys \cite{BhattBMS}, \cite{MorrowBMS} for details on this construction.

\section{Shtukas for $\mathrm{Spec}~\Z$}

Roughly, the upshot of the $A_\inf$-cohomology theory is that the ``universal'' $p$-adic cohomology theory is given by a shtuka relative to $\Spf \Z_p$. Currently, the author is trying to understand to what extent it might be true that the ``universal'' cohomology theory is given by a shtuka relative to $\Spec \Z$. It seems that this is a very fruitful philosophy.

Let us phrase the question more precisely. Given a proper smooth scheme $X$ over a base scheme $S$, what is ``the'' cohomology of $X$? Experience in arithmetic and algebraic geometry shows that there is no simple answer, and that in fact there are many cohomology theories: singular, de~Rham, \'etale, crystalline, etc.\footnote{We use the word ``cohomology theory'' in a loose sense, but it should be of Weil type and compatible with base change, so for example $H^1$ of an elliptic curve should be of rank $2$ over some ring, and base changing the elliptic curve should amount to a corresponding base change of its $H^1$ along an associated map of rings. In particular, ``absolute'' cohomologies such as motivic cohomology or syntomic cohomology are disregarded.} It helps to organize these cohomology theories according to two parameters: First, they may only defined for certain (geometric) points $s\in S$, and their coefficients may not be $\Z$ but only $\Z/n\Z$ or $\Z_\ell$ or $\R$.

This makes it natural to draw a picture of a space $``S\times \Spec \Z"$, where at a point $(s,\ell)$ we put $``H^i(X_s,\Z/\ell\Z)"$, whatever we want this to mean. Note that inside this space, we have the graph of $S$, taking any $s\in S$ to $(s,\mathrm{char}(s))$ where $\mathrm{char}(s)$ is the characteristic of $s$. For simplicity, we assume that $S=\Spec \Z$ in the following; then the graph is given by the diagonal $\Spec \Z\subset ``\Spec \Z\times \Spec \Z"$. We will use $p$ to denote a point of $p\in S=\Spec \Z$, and $\ell$ to denote a point of the vertical $\Spec \Z$.

Thus, let $X$ be a proper smooth scheme over $S=\Spec \Z$.\footnote{The known examples are not too interesting, so the reader may prefer to take $S=\Spec \Z[\tfrac 1N]$ for some $N$ to get more interesting examples. He should then remove the vertical fibres over $p|N$ from the picture.} We include the following cohomology theories in the picture:

\begin{figure}
\centering

\begin{tikzpicture}[auto]
\draw[->, very thick] (0,0) -- (0,9.5);
\draw[->, very thick] (0,0) -- (9.5,0);
\filldraw (0,10) circle [radius=.1];
\draw (-0.5,10) node {$\infty$};
\filldraw (10,0) circle [radius=.1];
\draw (10,-0.5) node {$\infty$};
\draw (8.5,-0.75) node {$\Spec \Z$};
\draw (5,-1) node {$p$};
\draw (-0.75, 8) node {$\Spec \Z$};
\draw (-1,5) node {$\ell$};
\filldraw (0,1.5) circle [radius=.1];
\draw (-0.5,1.5) node {$2$};
\filldraw (0,3) circle [radius=.1];
\draw (-0.5,3) node {$3$};
\filldraw (0,4.7) circle [radius=.1];
\draw (-0.5,4.7) node {$5$};
\filldraw (0,6.5) circle [radius=.1];
\draw (-0.5,6.5) node {$7$};
\filldraw (0,8) circle [radius=.05];
\filldraw (0,8.5) circle [radius=.05];
\filldraw (0,9) circle [radius=.05];
\filldraw (1.5,0) circle [radius=.1];
\draw (1.5,-0.5) node {$2$};
\filldraw (3,0) circle [radius=.1];
\draw (3,-0.5) node {$3$};
\filldraw (4.7,0) circle [radius=.1];
\draw (4.7,-0.5) node {$5$};
\filldraw (6.5,0) circle [radius=.1];
\draw (6.5,-0.5) node {$7$};
\filldraw (8,0) circle [radius=.05];
\filldraw (8.5,0) circle [radius=.05];
\filldraw (9,0) circle [radius=.05];
\draw[very thick, color=orange] (10,1) -- (10,10);
\draw[color=orange] (10.7,5) node {$\mathrm{singular}$};
\draw[very thick, color=green] (1,1) -- (10,10);
\draw[color=green] (8.9,8) node {$\mathrm{de~Rham}$};
\fill [blue, path fading=north] (1.8,1.5) rectangle (10,1.7);
\fill [blue, path fading=south] (1.8,1.3) rectangle (10,1.5);
\fill [blue, path fading=north] (1,3) rectangle (2.7,3.2);
\fill [blue, path fading=south] (1,2.8) rectangle (2.7,3);
\fill [blue, path fading=north] (3.3,3) rectangle (10,3.2);
\fill [blue, path fading=south] (3.3,2.8) rectangle (10,3);
\fill [blue, path fading=north] (1,4.7) rectangle (4.4,4.9);
\fill [blue, path fading=south] (1,4.5) rectangle (4.4,4.7);
\fill [blue, path fading=north] (5,4.7) rectangle (10,4.9);
\fill [blue, path fading=south] (5,4.5) rectangle (10,4.7);
\fill [blue, path fading=north] (1,6.5) rectangle (6.2,6.7);
\fill [blue, path fading=south] (1,6.3) rectangle (6.2,6.5);
\fill [blue, path fading=north] (6.8,6.5) rectangle (10,6.7);
\fill [blue, path fading=south] (6.8,6.3) rectangle (10,6.5);
\draw[color=blue] (8.3,5.5) node {$\mathrm{\acute{e}tale}$};
\draw[very thick, color=red] (1.5,1.3) -- (1.5,1.7);
\draw[very thick, color=red] (3,2.8) -- (3,3.2);
\draw[very thick, color=red] (4.7,4.5) -- (4.7,4.9);
\draw[very thick, color=red] (6.5,6.3) -- (6.5,6.7);
\draw[color=red] (5.2,4.1) node {$\mathrm{crystalline}$};
\end{tikzpicture}

\caption{A picture of known cohomology theories.}
\end{figure}
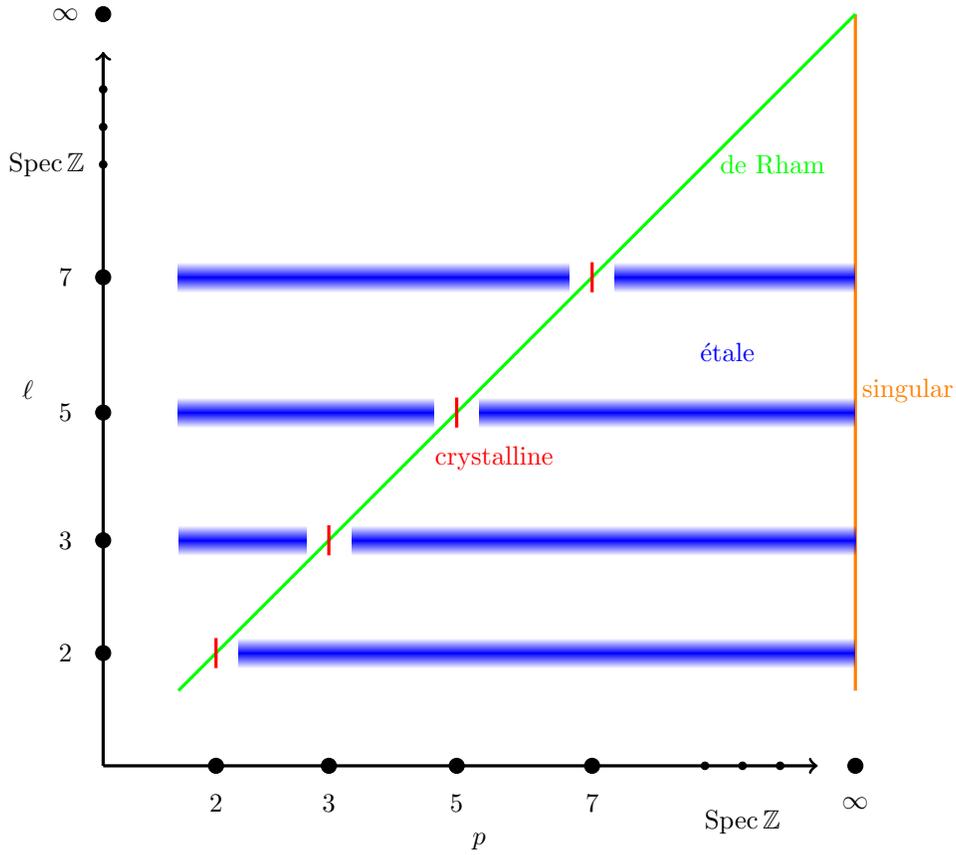

\begin{enumerate}
\item[{\rm (i)}] Singular cohomology $H^i_{\mathrm{sing}}(X(\C),\Z)$. This gives a vertical line over the infinite point of $\Spec \Z$.
\item[{\rm (ii)}] De~Rham cohomology $H^i_{\mathrm{dR}}(X/\Z)$. Note that (up to torsion problems)
\[
H^i_{\mathrm{dR}}(X/\Z)\otimes_{\Z} \Z/\ell\Z = H^i_{\mathrm{dR}}(X_{\F_p}/\F_p)
\]
for $p=\ell$, so a horizontal fibre of $H^i_{\mathrm{dR}}(X/\Z)$ agrees with a vertical fibre. Thus, it sits on the diagonal. It meets singular cohomology, via the comparison isomorphism
\[
H^i_{\mathrm{dR}}(X/\Z)\otimes_{\Z} \C\cong H^i_{\mathrm{sing}}(X(\C),\Z)\otimes_{\Z} \C\ .
\]
\item[{\rm (iii)}] \'Etale cohomology $H^i_\et(X_{\overline{s}},\Z/\ell^m\Z)$, which is defined for all geometric points $\overline{s}\in \Spec \Z\setminus \{\ell\}$. This gives a horizontal line with a hole at $\ell$, and in fact an infinitesimal neighborhood of the horizontal line. It meets singular cohomology via the comparison isomorphism
\[
H^i_{\et}(X_{\C},\Z/\ell^m\Z)\cong H^i_{\mathrm{sing}}(X(\C),\Z/\ell^m\Z)\ .
\]
\item[{\rm (iv)}] Crystalline cohomology $H^i_{\mathrm{crys}}(X_{\F_p}/\Z_p)$. This sits in a vertical fibre over $p$, and extends de~Rham cohomology infinitesimally.
\end{enumerate}

If we zoom in near a point $(p,p)$ in this picture, we arrive at a picture that looks exactly like Figure~\ref{fig:Ainf} depicting the $A_\inf$-cohomology theory! And indeed, one can use the $A_\inf$-cohomology theory to fill in this part of the picture.

The picture bears a remarkable similarity with the following equal-characteristic structure.

\begin{definition}[Anderson, Goss, \cite{Anderson}, \cite{Goss}] For a scheme $S=\Spec R$ over $\Spec \F_p[T]$ sending $T$ to $t\in R$, a t-motive is a finite projective $R[T]$-module $M$ together with an isomorphism
\[
\varphi_M: \Frob_R^\ast M[\tfrac 1{t-T}]\cong M[\tfrac 1{t-T}]\ .
\]
\end{definition}

Let us briefly discuss the similar specializations.

\begin{enumerate}
\item[{\rm (i)}] If $R=K$ is an algebraically closed field mapping to the infinite point of $\Spec \F_p[T]$ (which, strictly speaking, is not allowed -- but we can compactify $\Spec \F_p[T]$ into $\mathbb P^1_{\F_p}$) we get a vector bundle on $\P^1_K$ together with an isomorphism with its Frobenius pullback.\footnote{One might think that one has to allow a pole at $\infty$, but the analogy seems to suggest that it is in fact not there.} These are equivalent to vector bundles on $\P^1_{\F_p}$, and by restricting back to $\Spec \F_p[T]$, we get an $\F_p[T]$-module.
\item[{\rm (ii)}] Restricting along $R[T]\to R$ sending $T$ to $t$ we get a finite projective $R$-module. Moreover, $\varphi_M$ gives rise to a $R[[t-T]]$-lattice in $M\otimes_{R[T]} R((t-T))$, which gives rise to a filtration of $M\otimes_{R[T]} R$ similarly to the discussion in the case of $B_\dR^+$-lattices or $\C[[t]]$-lattices in Section~\ref{sec:twistor}.\footnote{A subtle point is that $M\otimes_{R[T],T\mapsto t} R$ is the analogue of Hodge(-Tate) cohomology, and de~Rham cohomology corresponds to the specialization $M\otimes_{R[T],T\mapsto t^p} R$.}
\item[{\rm (iii)}] Restricting modulo a power of $T$ (or $T-a$, $a\in \F_p$), the pair
\[
(M_n,\varphi_{M_n})=(M\otimes_{R[T]} R[t^{-1},T]/T^n,\varphi_M)
\]
is equivalent to an \'etale $\F_p[T]/T^n$-local system on $\Spec R[t^{-1}]$ by taking the \'etale sheaf $M_n^{\varphi_{M_n}=1}$.
\item[{\rm (iv)}] Restricting to the fibre $t=0$ (or $t=a$, $a\in \F_p$) and taking the $T$-adic completion, $M\otimes_{R[T]} (R/t)[[T]]$ with $\varphi_M$ defines a structure resembling crystalline cohomology.
\item[{\rm (v)}] Similarly, by the analogy between local shtukas and local $p$-adic shtukas, the picture near $(t,T)=(0,0)$ resembles the picture of the $A_\inf$-cohomology theory.
\end{enumerate}

An interesting question is what happens in vertical fibres, i.e.~for $S=\Spec \overline{\F}_p$. In the function field case, we get generically the following structure.

\begin{definition} A (generic) shtuka over $\Spec \overline{\F}_p$ relative to $\Spec \F_p(T)$ is a finite-dimensional $\overline{\F}_p(T)$-vector space $V$ together with a $\Frob_{\overline{\F}_p}$-semilinear isomorphism $\varphi_V: V\cong V$.
\end{definition}

Amazingly, Kottwitz \cite{KottwitzBG} was able to define an analogue of this category for number fields.

\begin{construction}[Kottwitz] For any local or global field $F$, there is an $F$-linear $\otimes$-category $\mathrm{Kt}_F$, constructed as the category of representations of a gerbe banded by an explicit (pro-)torus which is constructed using (local or global) class field theory. Moreover, one has the following identifications.
\begin{enumerate}
\item[{\rm (i)}] If $F=\F_p((T))$, then $\Kt_F$ is the category of finite-dimensional $\overline{\F}_p((T))$-vector spaces $V$ with a semilinear isomorphism $V\cong V$.
\item[{\rm (ii)}] If $F=\F_p(T)$, then $\Kt_F$ is the category of finite-dimensional $\overline{\F}_p(T)$-vector spaces $V$ with a semilinear isomorphism $V\cong V$.
\item[{\rm (iii)}] If $F=\Q_p$, then $\Kt_F$ is the category of finite-dimensional $W(\overline{\F}_p)[\tfrac 1p]$-vector spaces $V$ with a semilinear isomorphism $V\cong V$.
\item[{\rm (iv)}] If $F=\R$, then $\Kt_F$ is the category of finite-dimensional $\C$-vector spaces $V$ together with a grading $V=\bigoplus_{i\in \Z} V_i$ and a graded antiholomorphic isomorphism $\alpha: V\cong V$ such that $\alpha^2 = (-1)^i$ on $V_i$.
\end{enumerate}
\end{construction}

\begin{remark} In the local cases (i) and (iii), one can replace $\overline{\F}_p$ with any algebraically closed field of characteristic $p$ without changing the category. However, in the global case (ii), the equivalence holds only for $\overline{\F}_p$; for this reason, we fix this choice.
\end{remark}

It is an important problem to find a linear-algebraic description of $\Kt_\Q$. There are natural global-to-local maps, so there are functors $\Kt_\Q\to \Kt_{\Q_p}$ and $\Kt_\Q\to \Kt_\R$. The analogy between cohomology theories and shtukas suggests the following conjecture.

\begin{conjecture} There is a Weil cohomology theory $H^i_{\Kt_\Q}(X)$ for varieties $X$ over $\overline{\F}_p$ taking values in $\Kt_\Q$. Under the functor $\Kt_\Q\to \Kt_{\Q_p}$, this maps to crystalline cohomology, and under the functor $\Kt_\Q\to \Kt_{\Q_\ell}$ for $\ell\neq p$, this maps to \'etale cohomology (considered as an object of $\Kt_{\Q_\ell}$ via the fully faithful embedding of finite-dimensional $\Q_\ell$-vector spaces into $\Kt_{\Q_\ell}$). Under the restriction $\Kt_\Q\to \Kt_\R$, it gives a Weil cohomology theory $H^i_{\Kt_\R}(X)$ with values in $\Kt_\R$.
\end{conjecture}

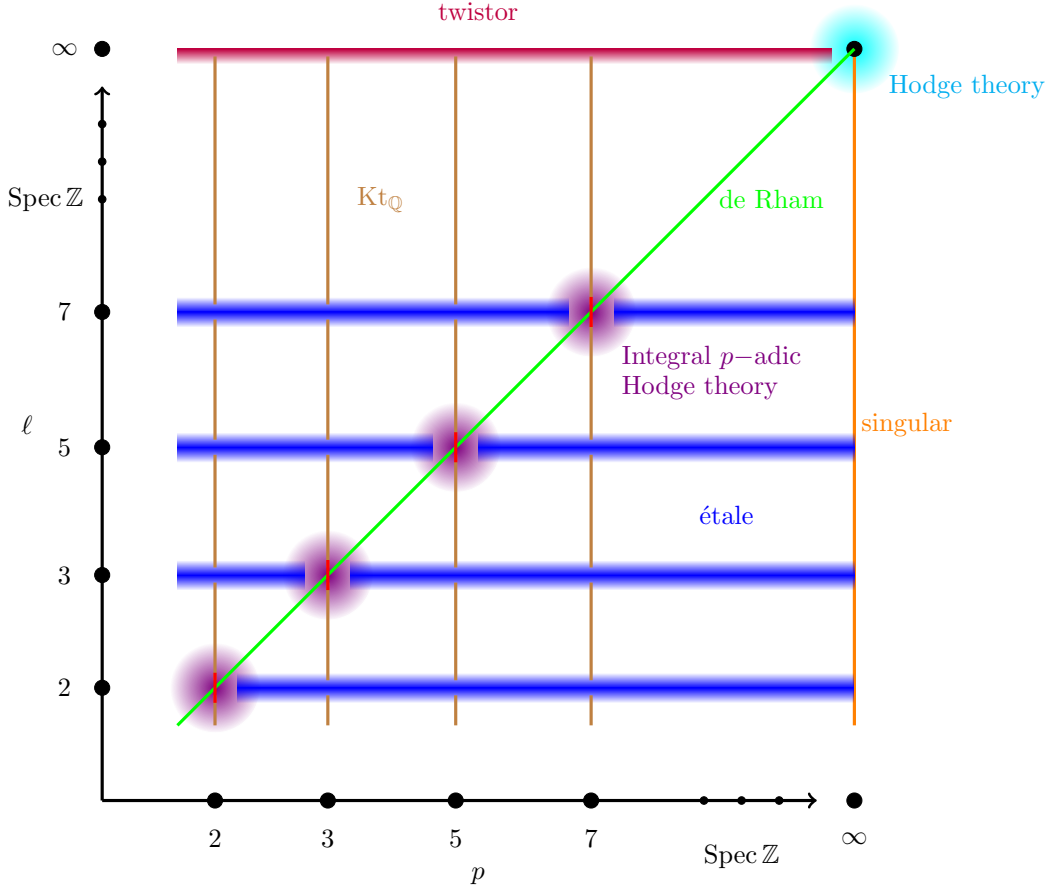
\begin{figure}
\centering

\begin{tikzpicture}[auto]
\draw[inner color=violet, color=white] (1.5,1.5) circle [radius=.6];
\draw[inner color=violet, color=white] (3,3) circle [radius=.6];
\draw[inner color=violet, color=white] (4.7,4.7) circle [radius=.6];
\draw[inner color=violet, color=white] (6.5,6.5) circle [radius=.6];
\draw[inner color=cyan, color=white] (10,10) circle [radius=.6];
\draw[color=cyan] (11.5,9.5) node {$\mathrm{Hodge~theory}$};
\draw[->, very thick] (0,0) -- (0,9.5);
\draw[->, very thick] (0,0) -- (9.5,0);
\filldraw (0,10) circle [radius=.1];
\draw (-0.5,10) node {$\infty$};
\filldraw (10,0) circle [radius=.1];
\draw (10,-0.5) node {$\infty$};
\draw (8.5,-0.75) node {$\Spec \Z$};
\draw (5,-1) node {$p$};
\draw (-0.75, 8) node {$\Spec \Z$};
\draw (-1,5) node {$\ell$};
\filldraw (0,1.5) circle [radius=.1];
\draw (-0.5,1.5) node {$2$};
\filldraw (0,3) circle [radius=.1];
\draw (-0.5,3) node {$3$};
\filldraw (0,4.7) circle [radius=.1];
\draw (-0.5,4.7) node {$5$};
\filldraw (0,6.5) circle [radius=.1];
\draw (-0.5,6.5) node {$7$};
\filldraw (0,8) circle [radius=.05];
\filldraw (0,8.5) circle [radius=.05];
\filldraw (0,9) circle [radius=.05];
\filldraw (1.5,0) circle [radius=.1];
\draw (1.5,-0.5) node {$2$};
\filldraw (3,0) circle [radius=.1];
\draw (3,-0.5) node {$3$};
\filldraw (4.7,0) circle [radius=.1];
\draw (4.7,-0.5) node {$5$};
\filldraw (6.5,0) circle [radius=.1];
\draw (6.5,-0.5) node {$7$};
\filldraw (8,0) circle [radius=.05];
\filldraw (8.5,0) circle [radius=.05];
\filldraw (9,0) circle [radius=.05];
\draw[very thick, color=orange] (10,1) -- (10,10);
\draw[color=orange] (10.7,5) node {$\mathrm{singular}$};
\filldraw (10,10) circle [radius=.1];
\draw[very thick, color=green] (1,1) -- (10,10);
\draw[color=green] (8.9,8) node {$\mathrm{de~Rham}$};
\fill [blue, path fading=north] (1.8,1.5) rectangle (10,1.7);
\fill [blue, path fading=south] (1.8,1.3) rectangle (10,1.5);
\fill [blue, path fading=north] (1,3) rectangle (2.7,3.2);
\fill [blue, path fading=south] (1,2.8) rectangle (2.7,3);
\fill [blue, path fading=north] (3.3,3) rectangle (10,3.2);
\fill [blue, path fading=south] (3.3,2.8) rectangle (10,3);
\fill [blue, path fading=north] (1,4.7) rectangle (4.4,4.9);
\fill [blue, path fading=south] (1,4.5) rectangle (4.4,4.7);
\fill [blue, path fading=north] (5,4.7) rectangle (10,4.9);
\fill [blue, path fading=south] (5,4.5) rectangle (10,4.7);
\fill [blue, path fading=north] (1,6.5) rectangle (6.2,6.7);
\fill [blue, path fading=south] (1,6.3) rectangle (6.2,6.5);
\fill [blue, path fading=north] (6.8,6.5) rectangle (10,6.7);
\fill [blue, path fading=south] (6.8,6.3) rectangle (10,6.5);
\draw[color=blue] (8.3,3.8) node {$\mathrm{\acute{e}tale}$};
\draw[color=violet] (8.1,5.9) node {$\mathrm{Integral~}p\mathrm{-adic}$};
\draw[color=violet] (7.95,5.5) node {$\mathrm{Hodge~theory}$};
\fill[purple, path fading=south] (1,9.8) rectangle  (9.7,10);
\draw[color=purple] (5,10.5) node {$\mathrm{twistor}$};
\draw[very thick, color=brown] (1.5,1) -- (1.5,1.4);
\draw[very thick, color=brown] (1.5,1.6) -- (1.5,2.9);
\draw[very thick, color=brown] (1.5,3.1) -- (1.5,4.6);
\draw[very thick, color=brown] (1.5,4.8) -- (1.5,6.4);
\draw[very thick, color=brown] (1.5,6.6) -- (1.5,9.9);
\draw[very thick, color=brown] (3,1) -- (3,1.4);
\draw[very thick, color=brown] (3,1.6) -- (3,2.9);
\draw[very thick, color=brown] (3,3.1) -- (3,4.6);
\draw[very thick, color=brown] (3,4.8) -- (3,6.4);
\draw[very thick, color=brown] (3,6.6) -- (3,9.9);
\draw[very thick, color=brown] (4.7,1) -- (4.7,1.4);
\draw[very thick, color=brown] (4.7,1.6) -- (4.7,2.9);
\draw[very thick, color=brown] (4.7,3.1) -- (4.7,4.6);
\draw[very thick, color=brown] (4.7,4.8) -- (4.7,6.4);
\draw[very thick, color=brown] (4.7,6.6) -- (4.7,9.9);
\draw[very thick, color=brown] (6.5,1) -- (6.5,1.4);
\draw[very thick, color=brown] (6.5,1.6) -- (6.5,2.9);
\draw[very thick, color=brown] (6.5,3.1) -- (6.5,4.6);
\draw[very thick, color=brown] (6.5,4.8) -- (6.5,6.4);
\draw[very thick, color=brown] (6.5,6.6) -- (6.5,9.9);
\draw[very thick, color=red] (1.5,1.3) -- (1.5,1.7);
\draw[very thick, color=red] (3,2.8) -- (3,3.2);
\draw[very thick, color=red] (4.7,4.5) -- (4.7,4.9);
\draw[very thick, color=red] (6.5,6.3) -- (6.5,6.7);
\draw[color=brown] (3.7,8) node {$\Kt_\Q$};
\end{tikzpicture}

\caption{A picture of the known, and some unknown, cohomology theories.}
\end{figure}

This conjecture is known to follow from the conjunction of the standard conjectures and the Tate conjecture (over $\overline{\F}_p$) and the Hodge conjecture (for CM abelian varieties) by the work of Langlands-Rapoport, \cite{LanglandsRapoport}. In fact, Langlands-Rapoport show that these conjectures imply that the category of motives over $\overline{\F}_p$ can be described as the category of representations of an explicit gerbe, and it is clear by inspection that there is a surjective map from Kottwitz' gerbe (which is independent of $p$!) to Langlands-Rapoport's gerbe. In particular, $\Kt_\Q$ conjecturally contains the category of motives over $\overline{\F}_p$ as a full subcategory, for all $p$.

It may be interesting to note how Serre's objection to a $2$-dimensional $\R$-cohomology theory is overcome by the $\R$-linear cohomology theory $H^i_{\Kt_\R}$, in an essentially minimal way: For a supersingular elliptic curve $E/\overline{\F}_p$, its associated $H^1_{\Kt_\R}(E)$ will be given by a $2$-dimensional $\C$-vector space $V=V_1$ together with an antiholomorphic isomorphism $\alpha: V\cong V$ such that $\alpha^2=-1$, so $\alpha$ does not give rise to a real structure; instead, it gives rise to a quaternionic structure. But the endomorphism algebra of $E$, which is a quaternion algebra over $\Q$ that is nonsplit over $\R$, can of course act on the Hamilton quaternions!

It is also interesting to note that there are functors from the category of isocrystals $\Kt_{\Q_p}$ to the category of vector bundles on the Fargues-Fontaine curve, and similarly from the category $\Kt_\R$ to the category of vector bundles on the twistor-$\P^1$, which in both cases induce a bijection on isomorphism classes. This gives another strong indication of the parallel between the Fargues-Fontaine curve and the twistor-$\P^1$, and how Kottwitz' categories play naturally into them.

In particular, one can also draw a horizontal line at $\ell=\infty$, where for all $s\in S$ one gets a vector bundle on the twistor-$\P^1$ (via the previous conjecture in finite characteristic, and via the twistor interpretation of Hodge theory in characteristic $0$). In other words, the complex variation of twistor structures must, in a suitable sense, be defined over the scheme $S$ (which is a general scheme over $\Spec \Z$, not necessarily over $\C$).

\section{$q$-de~Rham cohomology}

One may wonder what the completion of $``\Spec \Z\times \Spec \Z"$ along the diagonal looks like; in the picture of cohomology theories, this should combine $p$-adic Hodge theory for all primes $p$ with usual Hodge theory. One proposal for how this might happen was made in \cite{ScholzeqdeRham}. This paper builds on the following observation on the $A_\inf$-cohomology theory. Namely, one can write down explicit complexes computing this cohomology on affines. For example, if $R=\OO_C\langle T\rangle$ is the $p$-adic completion of $\OO_C[T]$, then one looks at the following complex:
\[
A_\inf\langle T\rangle \buildrel{\nabla_q}\over\longrightarrow A_\inf\langle T\rangle\ :\ T^n\mapsto \nabla_q(T^n)=[n]_q T^{n-1}\ ,
\]
where $[n]_q = 1+q+\ldots+q^{n-1} = \frac{q^n-1}{q-1}$ is Gau\ss 's $q$-integer, and $q=[\epsilon]\in A_\inf$ is a certain element of $A_\inf$. This is precisely the $q$-de~Rham complex studied by Aomoto, \cite{Aomoto}; it sends a general function $f(T)$ to the Jackson $q$-derivative
\[
(\nabla_q f)(T) = \frac{f(qT)-f(T)}{qT-T}\ .
\]
Unfortunately, the $q$-de~Rham complex depends heavily on the choice of coordinates, and it is not clear how to see the independence of the $A_\inf$-cohomology from the choice of coordinates using this picture.

However, in recent work with Bhatt, \cite{BhattScholzePrismatic}, we were able to prove the following theorem, proving a conjecture of \cite{ScholzeqdeRham}.

\begin{theorem}[\cite{BhattScholzePrismatic}]\label{thm:qdeRham} There is a $\Z[[q-1]]$-linear $q$-de~Rham cohomology theory for smooth schemes over $\Spec \Z$ that for any smooth $\Z$-algebra $R$ with a choice of an \'etale map $\Z[T_1,\ldots,T_n]\to R$ is computed by a $q$-deformation $q\Omega_R^\bullet$ of the de~Rham complex $\Omega^\bullet_R$. For example,
\[
q\Omega_{\Z[T]}^\bullet = \Z[T][[q-1]] \buildrel{\nabla_q}\over\longrightarrow \Z[T][[q-1]]\ : T^n \mapsto \nabla_q(T^n)=[n]_q T^{n-1}\ .
\]
After base change along $\Z[[q-1]]\to A_\inf$ via $q\mapsto [\epsilon]$, this recovers the $A_\inf$-cohomology theory.
\end{theorem}

In particular, this suggests that the completion of $``\Spec \Z\times \Spec \Z"$ along the diagonal is related to $\Spf \Z[[q-1]]$.

Previous progress towards this result was made by Pridham, \cite{Pridham}, and Cha\-tzistamatiou, \cite{Chatzistamatiou}, and announced by Masullo. In particular, Chatzistamatiou was able to write down explicitly the quasi-automorphism of $q\Omega_{\Z[T]}^\bullet$ for an automorphism of $\Z[T]$, like $T\mapsto T+1$; this is a highly nontrivial result!\footnote{It is virtually impossible to follow the constructions of \cite{BMS} to get an explicit quasi-isomorphism over $A_\inf$.}

Unfortunately, the proof of Theorem~\ref{thm:qdeRham} proceeds by first constructing such complexes after $p$-adic completion, and then patching them together in some slightly artificial way, so we believe that there is still much more to be understood about what exactly happens along the diagonal.

In fact, the prismatic cohomology defined in \cite{BhattScholzePrismatic} gives in the $p$-adic case a Breuil-Kisin variant\footnote{Taking values in Breuil-Kisin modules as defined in \cite{Kisin}.} of the $A_\inf$-cohomology theory for varieties over $\OO_K$ where $K$ is a finite extension of $\Q_p$. This Breuil-Kisin module contains some slightly finer information than the $q$-de~Rham cohomology and is in some sense a Frobenius descent of it, but it is not clear how to combine them for varying $p$. This is related to the observation that in the analogy with t-motives, the diagonal $T=t$ corresponds to Hodge(-Tate) cohomology, while the Frobenius-twisted diagonal $T=t^p$ corresponds to de~Rham cohomology. In this picture, the $q$-de~Rham cohomology would sit at the completion at $T=t^p$, while the prismatic picture seems more adapted to a picture corresponding to the completion at the diagonal $T=t$.

\bibliographystyle{amsalpha}
\bibliography{ICM}

\end{document}